\documentclass{svproc}
\usepackage{url}

\usepackage{hyperref}
\usepackage{type1cm}        
\usepackage{makeidx}         
\usepackage{graphicx}        

\usepackage{multicol}        
\usepackage[bottom]{footmisc}

\usepackage{newtxtext}       %
\usepackage[varvw]{newtxmath}       

\usepackage{bm}
\usepackage{siunitx}

 \usepackage{mathtools}

\newcounter{algorithm}

\newcommand{\bsa}{{\boldsymbol{a}}}

\newcommand{\bsc}{{\boldsymbol{c}}}
\newcommand{\bsd}{{\boldsymbol{d}}}

\newcommand{\bsh}{{\boldsymbol{h}}}

\newcommand{\bsk}{{\boldsymbol{k}}}

\newcommand{\bsn}{{\boldsymbol{n}}}

\newcommand{\bst}{{\boldsymbol{t}}}

\newcommand{\bsx}{{\boldsymbol{x}}}

\newcommand{\bsC}{{\boldsymbol{C}}}

\newcommand{\bsJ}{{\boldsymbol{J}}}

\newcommand{\bsT}{{\boldsymbol{T}}}

\newcommand{\bsV}{{\boldsymbol{V}}}

\newcommand{\bsX}{{\boldsymbol{X}}}

\newcommand{\bsell}{{\boldsymbol{\ell}}}
\newcommand{\bszero}{{\boldsymbol{0}}} 
\newcommand{\bsone}{{\boldsymbol{1}}}  

\newcommand{\bsgamma}{{\boldsymbol{\gamma}}}

\newcommand{\bsDelta}{{\boldsymbol{\Delta}}}

\newcommand{\bsPsi}{{\boldsymbol{\Psi}}}

\newcommand{\rmd}{{\mathrm{d}}}

\newcommand{\bbE}{{\mathbb{E}}}

\newcommand{\NN}{{\mathbb{N}}} 
\newcommand{\RR}{{\mathbb{R}}} 

\DeclareSymbolFont{bbold}{U}{bbold}{m}{n}
\DeclareSymbolFontAlphabet{\mathbbold}{bbold}
\newcommand{\ind}{{\mathbbold{1}}}

\newcommand{\calH}{{\mathcal{H}}}

\newcommand{\calK}{{\mathcal{K}}}

\newcommand{\calM}{{\mathcal{M}}}
\newcommand{\calN}{{\mathcal{N}}}
\newcommand{\calO}{{\mathcal{O}}}

\newcommand{\calU}{{\mathcal{U}}}

\DeclareMathOperator{\var}{Var}

\newcommand{\figpath}{Figures}
\usepackage{xspace}

\begin{document}

\mainmatter              
\title{Quasi-Monte Carlo Methods:  What, Why, and How?}
\titlerunning{Quasi-Monte Carlo Methods}  
%
\author{Fred J. Hickernell\inst{1} \and Nathan Kirk\inst{2} \and Aleksei G. Sorokin\inst{3}}
\authorrunning{Hickernell, Kirk, and Sorokin} 
%
\tocauthor{Fred J. Hickernell, Nathan Kirk, Aleksei G. Sorokin}
\institute{
Department of Applied Mathematics, Illinois Institute of Technology, \\ Chicago, IL, 60616, USA\\
\email{hickernell@iit.edu}, \
\texttt{http://www.iit.edu/~hickernell}
\and
\email{nkirk@iit.edu}
\and
\email{asorokin@hawk.iit.edu}
}

\maketitle              

\begin{abstract}
Many questions in  quantitative finance, uncertainty quantification, and other disciplines are answered by computing the population mean, $\mu := \mathbb{E}(Y)$, where instances of $Y:=f(\boldsymbol{X})$ may be generated by numerical simulation and $\bsX$ has a simple probability  distribution. The population mean can be approximated by the sample mean, $\hat{\mu}_n := n^{-1} \sum_{i=0}^{n-1} f(\bsx_i)$ for a well chosen sequence of nodes, $\{\bsx_0, \bsx_1, \ldots\}$ and a sufficiently large sample size, $n$.  Computing $\mu$ is equivalent to computing a $d$-dimensional integral, $\int f(\bsx) \varrho(\bsx) \, \mathrm{d} \bsx$, where $\varrho$ is the probability density for $\bsX$.

Quasi-Monte Carlo methods replace independent and identically distributed  sequences of random vector nodes, $\{\bsx_i \}_{i = 0}^{\infty}$, by low discrepancy sequences.  This accelerates the convergence of $\hat{\mu}_n$ to $\mu$ as $n \to \infty$.

This tutorial describes  low discrepancy sequences  and their quality measures.  We demonstrate the performance gains possible with quasi-Monte Carlo methods.  Moreover, we describe how to formulate problems to realize the greatest performance gains using quasi-Monte Carlo.  We also briefly describe the use of quasi-Monte Carlo methods for problems beyond computing the mean, $\mu$.

\end{abstract}
\setcounter{tocdepth}{2}

\section{Introduction} \label{sec:intro}
There are many settings where key underlying quantities that affect the outcome are unknown, e.g.,
\begin{itemize}
	\item Future market forces, which affect financial risk,
	\item The porosity field of rock, which affects the extraction of oil or gas, or
	\item Elementary particle interactions in a high energy physics experiment, which affect what is observed by detectors.
\end{itemize}
In such situations, the unknown inputs are often modeled using random vector variables or stochastic processes.  Computations are performed by generating a multitude of possible outcomes informed by the assumed probability distribution of the input. These are used to estimate the mean, quantiles, and/or probability distribution of the outcome.  This is the essence of the Monte Carlo (MC) method.

In mathematical terms, the random outcome is $Y := f(\bsX)$, where $\bsX$ is a vector random input.  Given an input, $\bsx$, the corresponding outcome  $y = f(\bsx)$ can be computed by an algorithm, whose complexity makes $f$ a \emph{black box}, or at least a \emph{gray box}.  The user selects a sequence of input values, also known as data sites or \emph{nodes}, $\bsx_0, \bsx_1, \ldots$, which lead to observed outcomes $y_0 = f(\bsx_0), y_1 = f(\bsx_1), \ldots$.  The $y_i$ may be used to approximate $\bbE(Y)$, quartiles of $Y$, the probability density of $Y$, and other quantities of interest.

Simple  Monte Carlo chooses the sequence of nodes, $\{\bsx_i\}_{i=0}^\infty$, to be independent and identically distributed (IID).  Quasi-Monte Carlo (qMC)  chooses the nodes differently, so that their empirical distribution approximates well the true probability distribution of $\bsX$.  The difference between  the true and empirical distributions is called a \emph{discrepancy}, and the node sequences used in qMC are called low discrepancy (LD) sequences.

This tutorial describes what qMC is, why we would want to use qMC, and how qMC can be implemented well.  The next section illustrates by example the advantages of qMC.  This is followed by a description of how deterministic LD sequences are constructed (Section \ref{sec:construct}) and randomizations that preserve their low discrepancy (Section \ref{sec:random}).  We describe various measures of discrepancy (Section \ref{sec:discrepancy}).  We explain how to decide what sample size, $n$, is sufficient to meet the user's error requirements (Section \ref{sec:stop}).  We then discuss how to rewrite the problem of interest in a qMC-friendly way (Section \ref{sec:reformulate}). We briefly discuss applications of qMC beyond computing the mean (Section \ref{sec:beyond}), and finish with a short conclusion (Section \ref{sec:conclusion}).

\section{An Illustration of Quasi-Monte Carlo} \label{sec:practice}

We illustrate the benefits of qMC with an example from Keister \cite{Kei96}, motivated by computational physics:

\begin{equation}\label{eq:keisterA}
	\mu := \int_{\mathbb{R}^d} \cos(\lVert \bst \rVert_2) \exp(-\lVert \bst \rVert_2^2) \, \rmd \bst,
\end{equation}
where $\lVert \bst \rVert_2 := \sqrt{t_1^2 + \cdots + t_d^2}$.  This integral may be evaluated numerically by re-writing it in spherical coordinates as
\begin{equation}\label{eq:keisterExact}
	\mu = \frac{2 \pi^{d/2}}{\Gamma(d/2)}\int_{0}^{\infty} \cos(r) \exp(-r^2) \, r^{d-1} \rmd r,
\end{equation}
where $2 \pi^{d/2}/\Gamma(d/2)$ is the surface area of the sphere in $d$ dimensions, and $\Gamma$ is the Gamma function.  The resulting one-dimensional integral is amenable to quadrature methods.  This non-trivial test case with a true value that can be easily calculated allows us to compute the numerical errors of various cubature schemes and contrast their performance.

For the purpose of this illustration, we work with $\mu$ in its original form, \eqref{eq:keisterA} rather than \eqref{eq:keisterExact}, and approximate it by a sample mean,
\begin{equation} \label{eq:sample_mean}
	\hat{\mu}_n := \frac 1n \sum_{i=0}^{n-1} f(\bsx_i).
\end{equation}
This does require some preparation to identify a suitable $f$, which is discussed in generality in Section \ref{sec:reformulate}.

This $\mu$ in \eqref{eq:keisterA} is the expectation of $Y := g(\bsT) := \pi^{d/2} \cos(\lVert \bsT \rVert)$, where $T_1, \ldots, T_d$ are IID Gaussian (normal) random variables with zero mean and variance $1/2$, i.e., $\bsT :=(T_1, \ldots, T_d) \sim \calN(\bszero,\mathsf{I}/2)$:
\begin{equation}\label{eq:keisterB}
	\mu = \int_{\mathbb{R}^d} \underbrace{\pi^{d/2} \cos(\lVert \bst \rVert)}_{=:g(\bst)}  \underbrace{\frac{\exp(-\lVert \bst \rVert^2)}{\pi^{d/2}}}_{\text{density of } \calN(\bszero,\mathsf{I}/2)} \, \rmd \bst = \bbE(Y) = \bbE[g(\bsT)].
\end{equation}
Virtually all LD sequences underlying qMC are defined to approximate the standard uniform distribution, $\calU[0,1]^d$.  Thus, we perform a variable transformation $\bst = \bigl (\Phi^{-1}(x_1), \ldots, \Phi^{-1}(x_d) \bigr )/\sqrt{2}$,
where is $\Phi$ is the cumulative distribution function of the standard Gaussian random variable.  This reimagines the integral $\mu$ as the expectation of a function, $f$, of a standard uniform random variable:
\begin{multline}\label{eq:keisterC}
	\mu = \int_{[0,1]^d} \underbrace{\pi^{d/2} \cos\Bigl (\bigl \lVert \bigl( \Phi^{-1}(x_1), \ldots, \Phi^{-1}(x_d)\bigr) \bigr\rVert/\sqrt{2}  \Bigr)}_{=:f(\bsx)} \, \rmd \bsx \\
	= \bbE[Y]
	= \bbE[f(\bsX)], \qquad \bsX \sim \calU[0,1]^d.
\end{multline}

\begin{figure}
	\centering
	\includegraphics[width =0.5\textwidth]{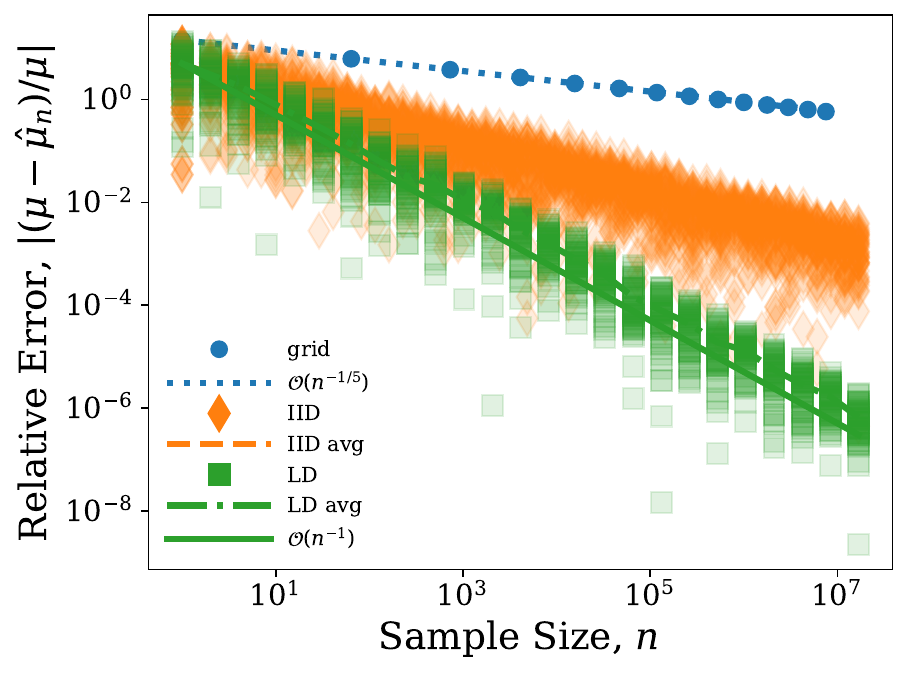}
	\caption{The relative error of approximating $\mu$ defined in \eqref{eq:keisterC} by the sample mean, $\hat{\mu}_n$, defined in \eqref{eq:sample_mean_Keister} for various choices of nodes, $\bsx_0, \bsx_1, \ldots$.  Grids have the largest error and LD nodes have the smallest error. \label{fig:keister-err}}
\end{figure}

Now we can apply IID MC, qMC, and other computational methods to approximate $\mu$ by the sample mean
\begin{equation} \label{eq:sample_mean_Keister}
	\hat{\mu}_n = \frac 1n \sum_{i=0}^{n-1} f (\bsx_i) = \frac 1n \sum_{i=0}^{n-1} \pi^{d/2} \cos\Bigl (\bigl \lVert \bigl( \Phi^{-1}(x_{i1}), \ldots, \Phi^{-1}(x_{id})\bigr) \bigr\rVert/\sqrt{2}  \Bigr).
\end{equation}
Consider the specific case of $d=6$ for which $\mu \approx -2.327303729298$.  We approximate $\mu$ by $\hat{\mu}_n$ for different, $\{\bsx_i\}_{i=0}^{n-1}$, for various sample sizes, $n$, and plot the relative errors, $\lvert (\mu - 	\hat{\mu}_n)/\mu\rvert$ in Figure \ref{fig:keister-err}. The three kinds of nodes are
\begin{enumerate}
	\renewcommand{\labelenumi}{\roman{enumi}.}
	\item Cartesian grids, $\{1/(2m), \ldots, (2m-1)/(2m) \}^6$, for $m = 2, 3, \ldots$ and $n = m^6$ (blue dots),
	\item IID sequences with arbitrary $n$ (orange diamonds), and
	\item Randomized LD (Sobol') sequences with $n = 1, 2, \ldots, 2^m, \ldots $ (green squares).
\end{enumerate}
There are $50$ replications for each of the latter two cases with the average errors being plotted as trend lines.

\begin{figure}
	\centering
	\includegraphics[width = \textwidth]{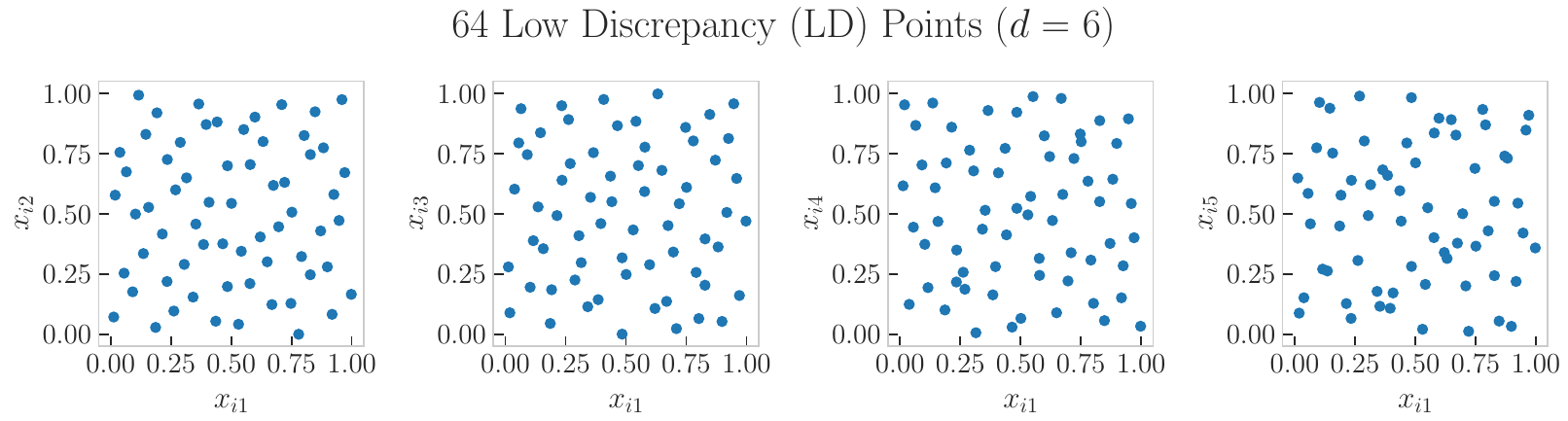}
	\caption{LD nodes cover the unit cube even better than IID nodes or grids.  In any  projection there is  a similar looking distribution of all $64$ nodes. \label{fig:ld}}
\end{figure}

Note the following from this example:
\begin{itemize}
	\item This integral is not particularly easy to evaluate numerically.  Even the best choice of nodes requires at least $n=100$ to get a relative error below $10\%$.

	\item \emph{Grids.} Although they may be attractive for low dimensional problems, e.g., $d = 1$, $2$, or $3$, grids do poorly for this modest dimension, $d=6$.  Figure \ref{fig:grid} compares a $64$ point grid with $d = 2$ and $6$.  For $d = 6$, the possible sample sizes, $n$, are quite sparse, as shown in Figure \ref{fig:keister-err}.

    Moreover, grids are not naturally extensible.  To move from $n = m^6$ nodes to $n = (m+1)^6$ nodes, one must discard the original nodes.

	Choosing the nodes to lie on a grid corresponds to a tensor product midpoint cubature rule, which would normally be expected to have an error of $\mathcal{O}(n^{-2/d})$  for general $d$.  The error decay of $\mathcal{O}(n^{-1/5})$ rather than $\mathcal{O}(n^{-2/6})$ for this example may be due to the lack of smoothness of $f$ near the boundaries of the unit cube.
\begin{figure}
	\centering
	\includegraphics[width = 0.3\textwidth]{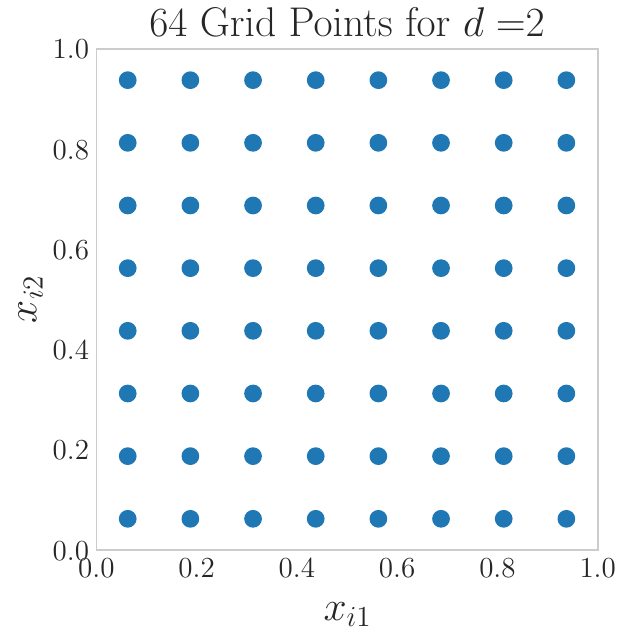}\qquad
	\includegraphics[width =  0.3\textwidth]{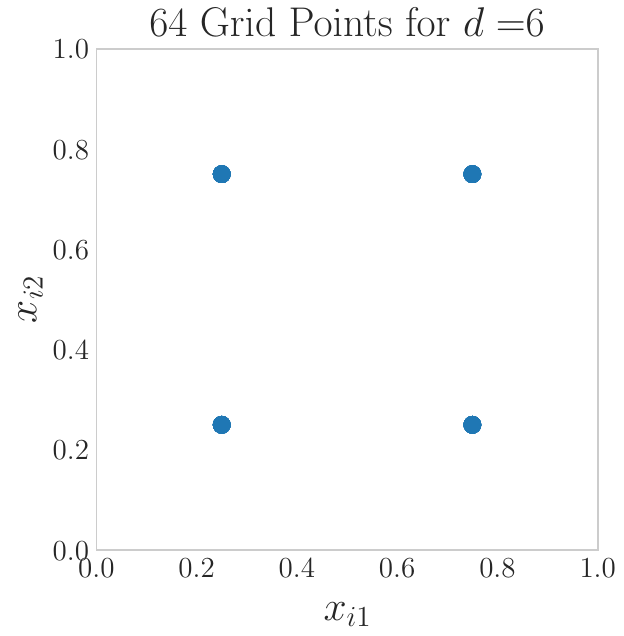}
	\caption{Although a two-dimensional grid covers the unit square rather well, a $d$-dimensional grid for modest $d$ does \emph{not} cover the unit cube well.  For example, one can only see four distinct nodes in a two dimensional projection of a six-dimensional grid with $64$ nodes. \label{fig:grid}}
\end{figure}

	\item \emph{IID.}  Simple MC is a substantial improvement over grid nodes. For IID nodes the root mean squared error is
	\begin{equation}\label{eq:IIDerror}
		\sqrt{\mathbb{E}[(\mu - \hat{\mu})^2]} = \sqrt{\var(\hat{\mu})} = \sqrt{\frac{\var(f(\bsX))}{n}} = \frac{\mathrm{Std}(f(\bsX))}{n^{1/2}},
	\end{equation}
	where $\var$ denotes the variance and $\mathrm{Std}$ the standard deviation.  This $\mathcal{O}(n^{-1/2})$ decay is observed in Figure \ref{fig:keister-err}. For simple MC the sample size, $n$, can be any positive integer without affecting the rate of decay of the error.

	Whereas grid points collapse on top of one another when viewed in low dimensional projections (Figure \ref{fig:grid}), all IID nodes may be seen when viewed in any lower dimensional projection, as seen in Figure \ref{fig:iid}.  The disadvantage of IID nodes is that they form clusters and leave gaps.  This is because the position of any one node is independent of the position of the others.

\begin{figure}
	\centering
	\includegraphics[width =\textwidth]{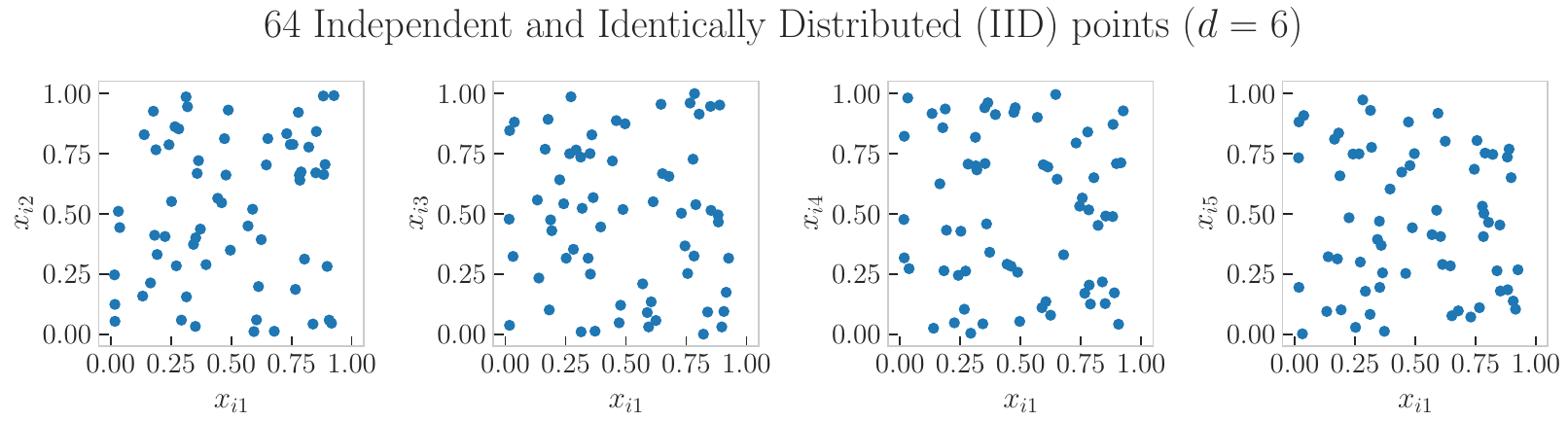}
	\caption{IID nodes cover the unit cube better than grid nodes, although one does observe clusters and gaps.  In any  projection there is  a similar looking distribution of all $64$ nodes. \label{fig:iid}}
\end{figure}

	\item \emph{LD.}  The error of qMC methods decays nearly like $\mathcal{O}(n^{-1})$, which for this example correponds to a reduction in error of several orders of magnitude compared to simple MC for large enough $n$.  Typically, LD sequences have preferred sample sizes.  This particular LD sequence is a Sobol' sequence, where the preferred sample sizes are non-negative integer powers of $2$.

	Figure \ref{fig:ld} shows typical two-dimensional projections of a size $64$ LD node set.  Visually, these nodes fill the unit cube better than IID nodes.  A quantitative measure of this, the discrepancy, is defined in Section \ref{sec:discrepancy}.  Constructions of LD node sequences are explained in the next section.

\end{itemize}

The numerical computations used to generate the figures in this chapter may be reproduced using the Jupyter notebook in \cite{MCQMC2024TutorialNotebook}. The Python library \texttt{qmcpy} \cite{QMCPy2020a} is used to generate the LD node sequences. 

\section{LD Sequence Constructions} \label{sec:construct}
This section introduces some of the most popular LD sequences.  What sets these apart from IID sequences is that the nodes are deliberately chosen and highly correlated, whether they be deterministic or randomized. The Python library \texttt{qmcpy} \cite{QMCPy2020a} contains many of the constructions discussed here, as well as the stopping criteria discussed in Section \ref{sec:stop}.

\subsection{Lattice Sequences} \label{sec:lattice}
\begin{figure}
	\centering
	\includegraphics[width = \textwidth]{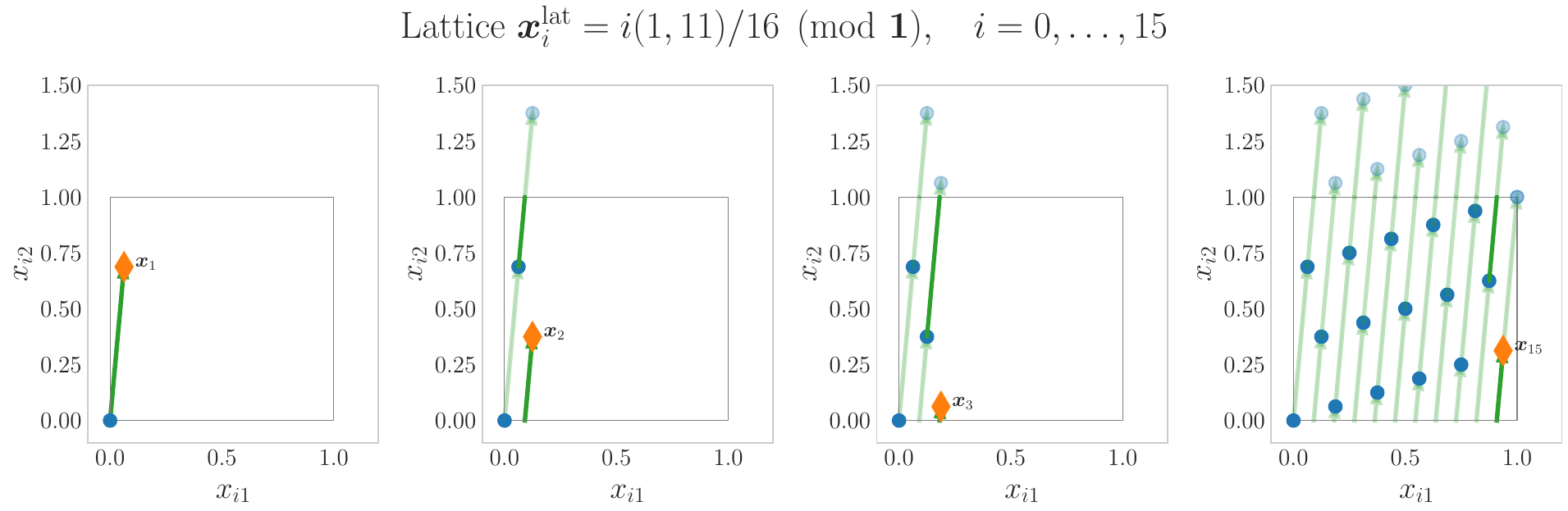}
	\caption{The construction of a good lattice node set in two dimensions with $16$ nodes is obtained by moving $\bsh/n$ beyond the present node and wrapping around the boundaries of the square (hypercube in $d$ dimensions) until one returns to the origin. \label{fig:latticeconstruct}}
\end{figure}

One of the simplest LD constructions is the family of good lattice points \cite{DicEtal22a,SloJoe94}.  As a finite node set, they are defined as
\begin{equation} \label{eq:latticepts}
	\bsx_i^{\text{lat}} = i\bsh/n \pmod \bsone, \qquad i = 0, \ldots, n-1,
\end{equation}
where $\bsh$ is a well chosen $d$-dimensional integer vector. Figure \ref{fig:latticeconstruct} illustrates this construction for $n = 16$ and $\bsh = (1,11)$.  Note that the set $\{\bsx^{\text{lat}}_0, \ldots, \bsx^{\text{lat}}_{n-1}\}$ defined in \eqref{eq:latticepts} is closed under addition modulo $\bsone$ so that it forms a group.

One disadvantage of this construction is that it is not extensible. For the example in Figure \ref{fig:latticeconstruct}, the first $8$ nodes do not fill the unit square well.  However, the set of nodes defined by even $i$ do a reasonable job.  This suggests a method for defining extensible lattice sequences that was proposed independently in \cite{Mai81a} and \cite{HicEtal00}.

An extensible lattice sequence in one dimension is defined by the van der Corput sequence in base $b$, $\{\phi_b(0), \phi_b(1), \ldots\}$. This involves the so-called radical inverse function and essentially reflects the digits of the integer $i$ in base $b$ about the decimal point. For example, $\phi_2(6) = \phi_2(110_{2}) = {}_20.011 = 3/8$.  In general,
\begin{multline} \label{eq:vdc}
	\phi_b(i_0 + i_1b + i_2 b^2 + \cdots ) := i_0 b^{-1} + i_1 b^{-2} + i_2 b^{-3} + \cdots \in [0,1)
	\\
	 \text{where } i_0, i_1, \ldots \in \{0,\ldots b-1\}.
\end{multline}
For all non-negative integers, $m$, the first $n = b^m$ nodes in the van der Corput sequence correspond to the evenly spaced nodes $\{0, b^{-m}, \ldots, 1 - b^{-m} \}$---albeit in a different order.

To construct an extensible lattice, we replace $i/n$ in \eqref{eq:latticepts} by $\phi_b(i)$ to get
\begin{equation} \label{eq:extensiblelattice}
	\bsx_i^{\text{lat}} =\phi_b( i)\bsh \pmod \bsone, \qquad i = 0, 1 , \ldots.
\end{equation}
This reordering of nodes from the original construction allows us to preserve the lattice structure for the first $b^m$ nodes for any non-negative integer $m$. That is, $\{\bsx^{\text{lat}}_0, \ldots, \bsx^{\text{lat}}_{b^m-1}\}$ is a closed under addition modulo $\bsone$.

\begin{figure}
	\centering
	\includegraphics[width = \textwidth]{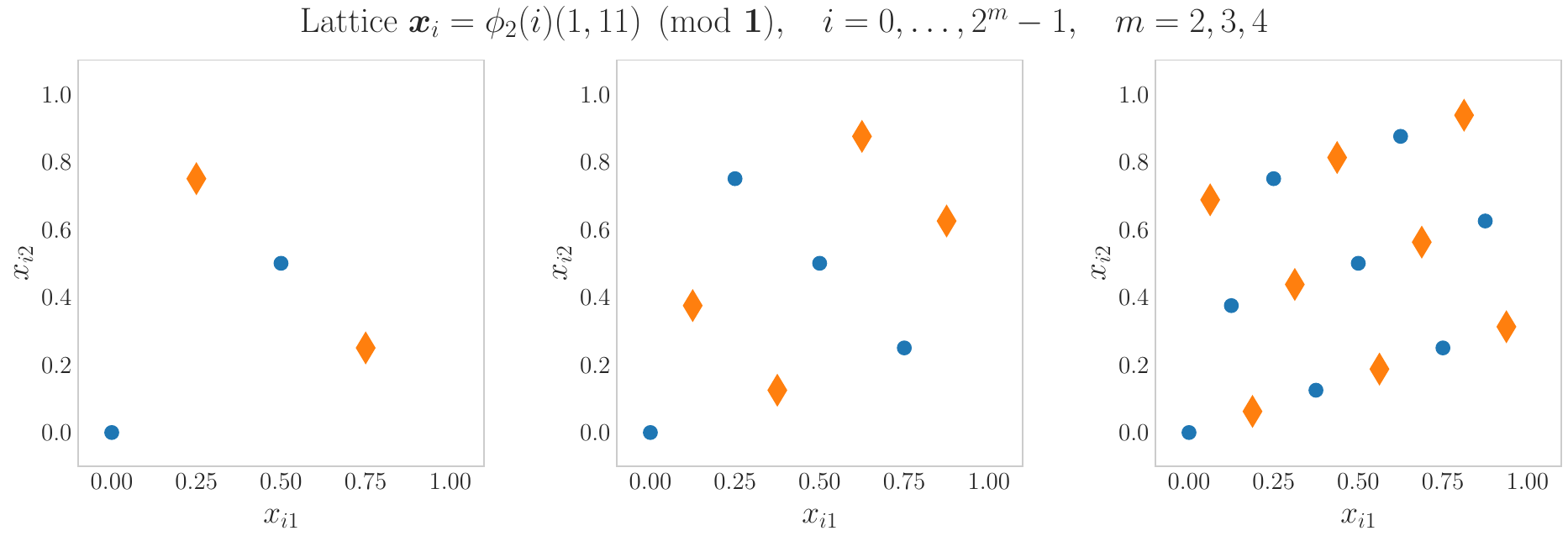}
	\caption{An extensible lattice corresponding to  Figure \ref{fig:latticeconstruct}  with the nodes reordered using the van der Corput sequence in base $2$.  For each plot the blue dots are a copy of the nodes to the left and the orange diamonds are a shifted copy (modulo $\bsone$) of the blue dots. \label{fig:extensiblelatticeconstruct}}
\end{figure}

Figure \ref{fig:extensiblelatticeconstruct} shows the first $4$, $8$, and $16$ nodes of an extensible lattice in base $2$ with the same generator as in Figure \ref{fig:latticeconstruct}.  The blue dots are a copy of the nodes in the plot to the left, and the the orange diamonds correspond to a shifted copy of the blue dots modulo $\bsone$.  For the left plot, the shift is $\phi_2(2)\bsh = (1,11)/4=(0.25,0.75)$.  For the middle plot, the shift is $\phi_2(4)\bsh = (1,11)/8=(0.125,0.375)$. For the right plot, the shift is $\phi_2(8)\bsh = (1,11)/16=(0.0625,0.6875)$.

For the example in Figure \ref{fig:extensiblelatticeconstruct}, increasing $n$ beyond $16$ repeats the original nodes because the generating vector, $\bsh$, contains integers all less than $16$.  Obtaining a truly extensible lattice sequence requires that $\bsh$ be a vector of generalized integers---essentially integers with infinite numbers of nonzero digits---as explained in \cite{HicNie03a}, where the existence of good generating vectors for infinite lattice sequences is also proved.  In practice, one searches computationally for a $\bsh$ that produces good set of lattice nodes of size $n = b^m$ for a range of $m$ \cite{HicEtal00}. The construction of great generating vectors for lattices has attracted a great deal of interest and is discussed further in Section \ref{sec:optLD}. 


\subsection{Digital Sequences} \label{sec:digital}

Another family of LD sequences that is a generalization of the van der Corput sequence \eqref{eq:vdc} is digital sequences \cite{DicPil10a,Nie92}.  For simplicity, we restrict ourselves to $b = 2$ and let $\oplus$ denote binary digitwise addition modulo $2$ (also known as digitwise exclusive or), e.g., $3/8 \oplus 3/4 = {}_20.011 \oplus {}_20.110 = {}_20.101 = 5/8$.  A digital sequence is defined as
\begin{equation} \label{eq:digital}
	\bsx_i^{\text{dig}} := i_0 \bsx^{\text{dig}} _1 \oplus i_1 \bsx^{\text{dig}}_2 \oplus i_2 \bsx^{\text{dig}} _4 \oplus \cdots \in [0,1)^d \quad \text{for }
	i = i_0 + i_12 + i_2 2^2 + \cdots,
\end{equation}
where $\bsx^{\text{dig}}_1, \bsx^{\text{dig}}_2, \ldots$ are carefully chosen.  (For $d=1$ and $x^{\text{dig}}_{2^m} = 2^{-m-1}$, this is the van  der Corput sequence in base $2$.) The node set $\{ \bsx^{\text{dig}}_0, \ldots, \bsx^{\text{dig}}_{2^m-1}\}$  for non-negative integer $m$ is called a digital net and is closed under $\oplus$.

If $x^{\text{dig}}_{ij\ell}$ denotes the $\ell^{\text{th}}$ digit of the $j^{\text{th}}$ coordinate of $\bsx^{\text{dig}}_i$, then the literature often defines generating matrices, $\mathsf{C}_j = \bigl(x^{\text{dig}}_{ij\ell}\bigr)_{\ell,i = 1}^{M,N}$ for $j = 1, \ldots, d$, where $M$ is the maximum number of bits in the expression for $\bsx^{\text{dig}}_i$, say $52$, and $2^N$ is the maximum number of nodes that is intended to be generated.  The digital sequence can then be defined equivalently as
\begin{equation} \label{eq:digitalB}
\renewcommand{\arraystretch}{1.5}
	\begin{pmatrix} x^{\text{dig}}_{ij1} \\ x^{\text{dig}}_{ij2} \\ \vdots \end{pmatrix}
	= \mathsf{C}_j \begin{pmatrix} i_0 \\ i_1 \\ \vdots \end{pmatrix}
    \quad \pmod{2}, \qquad j = 1,\ldots, d, \ i = 0, 1, \ldots.
\end{equation}

To understand how well these nodes are evenly distributed over $[0,1)^d$, imagine a box of the form
\begin{multline} \label{eq:element_box}
	[a_1 b^{-k_1}, (a_1 + 1) 2^{-k_1}) \times \cdots \times [a_d b^{-k_d}, (a_d + 1) 2^{-k_d}), \\
	a_j \in \{0, \ldots, 2^{k_j} - 1\}, \ k_j \in \{0,1, \ldots \}, \qquad j =1, \ldots d, \ \bsk \in \NN_0.
\end{multline}
This box has volume $2^{-(k_1 + \cdots +k_d)} = 2^{-\lVert\bsk\rVert_1}$.  For a digital net with $2^m$ nodes and $\lVert \bsk \rVert_1 \le m$, a ``fair share" of nodes for this box would be $2^{m- \lVert \bsk \rVert_1}$.  This will occur if and only if
\begin{multline} \label{eq:tcond}
	\text{the first $k_1$ rows of the first $m$ columns of $\mathsf{C}_1$ plus} \\
	\text{the first $k_2$ rows of the first $m$ columns of $\mathsf{C}_2$ plus} \\
	 \vdots \\
	\text{the first $k_d$ rows of the first $m$ columns of $\mathsf{C}_d$} \\
	\text{are linearly independent over addition and multiplication modulo $2$, }
\end{multline}
regardless of the choice of $\bsa = (a_1, \ldots, a_d)$.  The $t$-value of a digital net with $2^m$ nodes is defined as the smallest $t$ for which condition \eqref{eq:tcond} holds for all $\bsk$ with $\lVert \bsk \rVert_1 \le m-t$.  Equivalently, this $t$ is the smallest value for which every box of the form \eqref{eq:element_box} with volume $2^{\lVert \bsk \rVert_1 - m}$ contains its fair share of $2^t$ nodes.  Such a digital net is then called a $(t,m,d)$-net.  An infinite sequence of the form \eqref{eq:digitalB} for which this condition holds for all non-negative $m$ is called a  $(t,d)$-sequence.

\begin{figure}
	\centering
	\includegraphics[width = \textwidth]{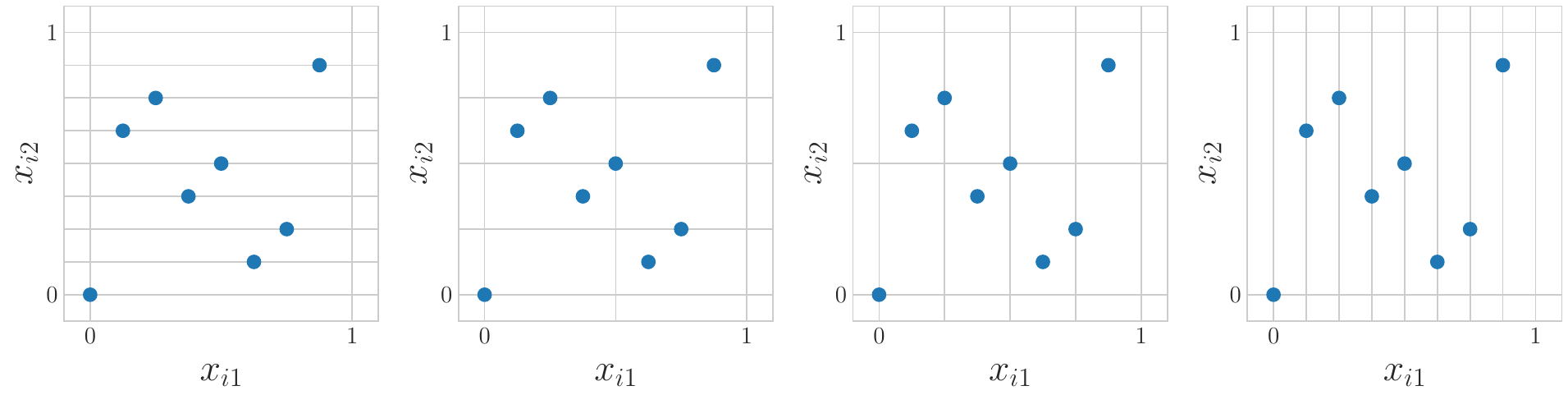}
	\includegraphics[width = \textwidth]{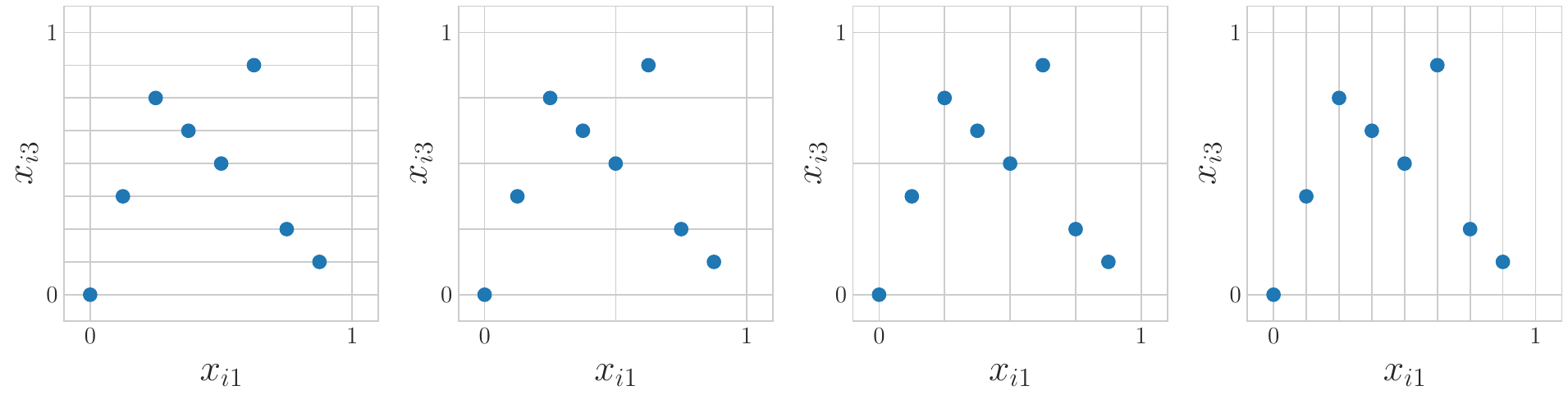}
	\caption{Each box in the two dimensional projections of the node set \eqref{eq:smallnet} plotted above contain one node except for the second row, third plot from the left.  Thus, this node set is not a $(0,3,3)$-net.  It is however a $(1,3,3)$-net.  \label{fig:elementinterval}}
\end{figure}

To illustrate the $t$-value, consider the following $d=3$ dimensional digital net with $2^3$ (eight) nodes, i.e., $m = 3$:
\begin{multline} \label{eq:smallnet}
	\{(0, 0,   0),
	(0.5,   0.5,   0.5  ),
	(0.25,  0.75,  0.75 ),
	(0.75,  0.25,  0.25 ),
	(0.125, 0.625, 0.375),\\
	(0.625, 0.125, 0.875),
	(0.375, 0.375, 0.625),
	(0.875, 0.875, 0.125)\}.
\end{multline}
Figure \ref{fig:elementinterval} shows several two dimensional projections of these nodes and the two-dimensional boxes (rectangles) of the form \eqref{eq:element_box} with $\lVert\bsk\rVert_1 = 3$.  Most of the boxes contain their fair share of one node, but the box $[0,1/4) \times [0,1) \times [0,1/2)$ contains two nodes and the adjacent box, $[1/4,1/2) \times [0,1) \times [0,1/2)$ contains no nodes.  Thus, the $t$-value cannot be $0$.  However, the $t$-value is $1$ because all boxes of the form  \eqref{eq:element_box} with $\lVert\bsk\rVert_1 = 2$, i.e., a volume of $2^{-2} = 1/4$ contain a fair share of $2^{m- \lVert \bsk \rVert_1} = 2$ nodes.  The node set \eqref{eq:smallnet} is a $(1,3,3)$-net.

This conclusion can also be reached by looking at the first three rows and columns of the generating matrices for this digital net, which are
\begin{equation*}
	\mathsf{C}_1 = \begin{pmatrix}
		1 & 0 & 0 & \cdots \\ 0 & 1 & 0 & \cdots \\ 0 & 0 & 1 & \cdots \\ \vdots & \vdots & \vdots & \ddots
	\end{pmatrix}, \quad
	\mathsf{C}_2 = \begin{pmatrix}
		1 & 1 & 1 & \cdots \\ 0 & 1 & 0 & \cdots \\ 0 & 0 & 1 & \cdots \\ \vdots & \vdots & \vdots & \ddots
	\end{pmatrix}, \quad
	\mathsf{C}_3 = \begin{pmatrix}
	1 & 1 & 0 & \cdots \\ 0 & 1 & 1 & \cdots \\ 0 & 0 & 1 & \cdots \\ \vdots & \vdots & \vdots & \ddots
\end{pmatrix}.
\end{equation*}
For $\bsk = (2, 0, 1)$ and $\lVert \bsk \rVert_1 = 3$, condition \eqref{eq:tcond} is not satisfied, and so $t$ cannot be $0$ for the node set  \eqref{eq:smallnet}.  However, condition \eqref{eq:tcond} is satisfied for all $\bsk$ with $\lVert \bsk \rVert_1 = 2$, again confirming that we have a $(1,3,3)$-net.

If we consider only the first two coordinates of the node set defined in \eqref{eq:smallnet}, then  condition \eqref{eq:tcond} is satisfied for all $\bsk$ with $\lVert \bsk \rVert_1 = 3$.  The first two coordinates of \eqref{eq:smallnet} are a $(0,3,2)$-net

Digital sequence generators can be found via number theory \cite[Chapter 8]{DicPil10a} or numerical search (see \cite{KuoJoe08a} and \cite[Chapter 10]{DicPil10a}).  The earliest instance is due to Sobol' \cite{Sob67}. 

\subsection{Halton sequences} \label{sec:Halton}  
\begin{figure}
	\centering
	\includegraphics[width=\textwidth]{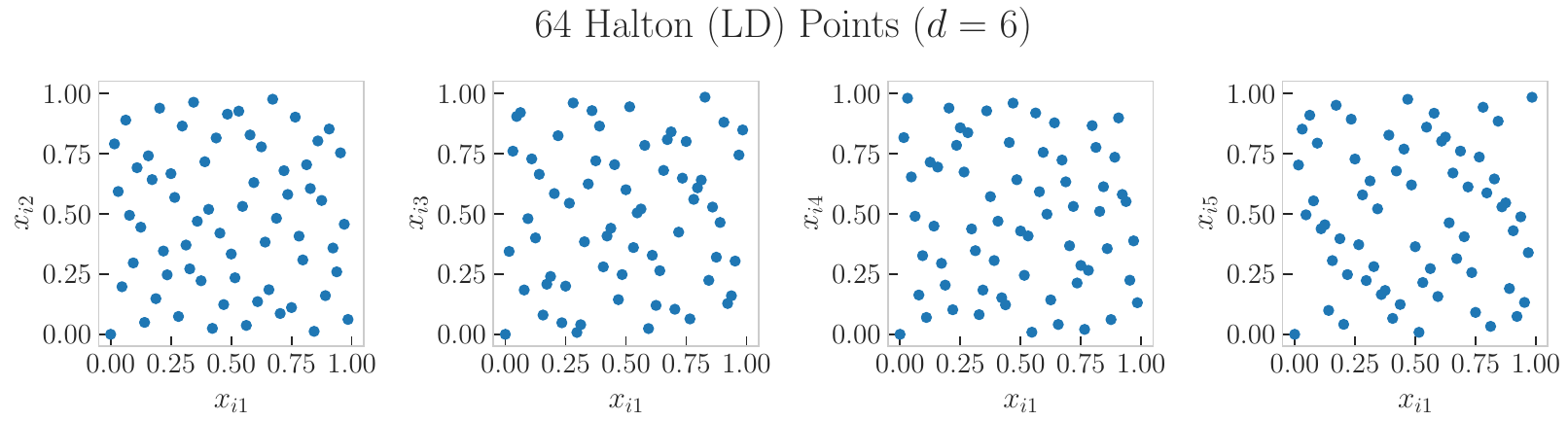}
	\caption{Halton nodes as defined in \eqref{eq:Halton} have low discrepancy but no preferred sample size. \label{fig:Halton}}
\end{figure}

While lattice sequences and digital sequences have preferred sample sizes, $n$, Halton sequences have no preferred sample size.  The Halton sequence is defined in terms of the van der Corput sequences for different bases:
\begin{equation}\label{eq:Halton}
	\bsx^{\text{Hal}}_i = \bigl( \phi_{b_1}(i), \ldots, \phi_{b_d}(i) \bigr), \qquad i = 0, 1, \ldots,
\end{equation}
where $b_1, \ldots, b_d$ is a choice of $d$ distinct prime numbers.  Often they are chosen as the first $d$ prime numbers.  Figure \ref{fig:Halton} shows two dimensional projections of Halton nodes in six dimensions. 

The Halton construction is an extensible LD sequence. However, should one know $n \in \mathbb{N}$ in advance, a finite $d+1$-dimensional LD node set called the Hammersley node set can be defined similarly as
\begin{equation}\label{eq:hammersley}
    \bsx_i^{\text{Ham}} = \left( \frac{i}{n}, \phi_{b_1}(i), \ldots, \phi_{b_d}(i) \right) \quad i = 0, \ldots, n-1.
\end{equation}

\subsection{LD Nodes by Optimization}\label{sec:optLD}

Most quasi-Monte Carlo constructions are studied based on their asymptotic behavior as \(n \to \infty\), with numerous results on asymptotic rates \cite{Nie92,NovWoz10a}. However, few methods address finding optimal node configurations for specific \(n\) and \(d\) which can be critical in practical applications where each function evaluation is quite costly. Of course, one possible strategy is to truncate asymptotically optimal sequences to the desired \(n\), however, as discussed above many of the standard LD constructions have a preferred number of nodes with truncations leading to imbalances in equidistribution. We can often do better optimizing directly for the target \(n\) and \(d\).

When optimizing LD nodes, one requires an objective function which in some sense must assess the quality of the LD nodes directly, or the expected performance when the node set is implemented in a qMC estimate \eqref{eq:sample_mean}.  We will refer to such an objective function as a figure-of-merit (FOM) and popular choices include the discrepancy (see Section \ref{sec:discrepancy}) or worst-case error metrics when restricted to a specific space of functions; see \cite{Hic00a,LatNet}.

For several decades, combinatorial optimization techniques have been employed to obtain good generating parameters for LD constructions. For example, in the case of a lattice rule which is fully determined by the generating vector \(\boldsymbol{h} \in \mathbb{Z}_n^d \coloneqq \{0, 1, \ldots, n-1\}^d\) given $n$ and $d$,
one often resorts to finding good $\boldsymbol{h}$ via computer search to minimize a chosen FOM. An exhaustive search is infeasible due to the exponential growth of the search space with \(d\). 

Thus, the component-by-component (CBC) construction, introduced by Korobov \cite{kor63} and later revisited by Sloan and collaborators \cite{Slo02}, reduces the search space to size \(dn\) by constructing \(\boldsymbol{h}\) one component at a time while keeping the previous components fixed. This greedy procedure also applies to finding good generating matrices $\mathsf{C}_j$ for digital sequences.
For certain choices of function spaces, there exist fast component-by-component construction implementations using fast Fourier transforms (FFTs) as the primary tool. We refer to \cite{DicKuo04a,JoeKuo03,KuoJoe02b,DicEtal22a} and references therein for works on the CBC construction for lattice rules and digital sequences, \cite{NuyCoo06a,NuyCoo06b} for their fast implementations and \cite{LEcEtal22a,LatNet} for a qMC software tool \texttt{LatNet Builder} by L'Ecuyer to generate LD nodes which relies heavily on the CBC construction method.

When constructing LD nodes, it is essential to ensure not only a well-distributed set of nodes in the $d$-dimensional space but also in its lower-dimensional projections (recall the Cartesian grid from Section \ref{sec:intro}, Figure \ref{fig:grid}). \texttt{MatBuilder} \cite{paulin2022} is a software tool developed precisely for this task; to generate generating matrices for digital sequence LD nodes possessing excellent low dimensional uniformity. In this software, the solutions to integer linear programming problems are used to build the columns of generator matrices in a greedy manner for the LD digital sequence.


We have outlined various efforts by qMC researchers to optimize specific classes of LD constructions, including lattices and digital sequences. Significant effort has even been taken to optimize the uniformity of so-called generalized Halton sequences due to the poor lower dimensional projections exhibited by even a moderate dimensional Halton sequence \cite{kirklem24}. These methods are discussed briefly in Section \ref{sec:randHalton}. There also exist several efforts to generate novel LD constructions which do not adhere to any underlying number theoretic structure, but none-the-less, possess a high level of uniformity. These methods are often global optimization methods of sample point distributions, which is often challenging due to the non-convex nature of the FOM as a function of the nodes. 
One of the first works to globally optimize LD nodes, albeit a heuristic approach, was by Winker and Fang \cite{WinFan97b} employing the threshold accepting optimization algorithm to directly optimize the uniformity of the LD nodes with respect to a discrepancy FOM. 

More recently, Rusch, Kirk, Bronstein, Lemieux and Rus \cite{ruschkirk24} developed the first machine-learning assisted optimization tool for LD node optimization in the Message-Passing Monte Carlo (MPMC) method. MPMC implements tools from geometric deep learning, imposing a computational graph onto an initialized IID node set to transform the nodes to LD node set via a learned mapping. The optimization is gradient-based and is guided by the direct minimization of a discrepancy-based FOM. 

We also highlight the method proposed by Hinrichs and Oettershagen \cite{hinoet16} generating optimal LD nodes for a specific class of functions in two-dimensions. To address the non-convexity challenge, the authors decompose the global optimization of a specific FOM into a large number of smaller, convex sub-problems which are solved individually. 

There are several other recent attempts to construct optimal LD nodes, many of which are well summarized in the PhD thesis of Cl\'{e}ment \cite{clethesis24}. We highlight the subset selection method proposed in \cite{cle22} to choose from an $n$-element initialization set, the $k<n$ nodes which minimize a discrepancy FOM. A heuristic approach to this problem was later shown in \cite{cle24_heuristic}. Furthermore, a method to generate optimal star-discrepancy node sets for fixed $n$ and $d$ based on a non-linear programming approach was suggested in \cite{cle24}, however is computationally restricted to very small numbers of nodes in dimensions two and three only.

\section{Randomization} \label{sec:random}
The LD constructions described in the previous section have so far been \emph{deterministic}. The sample mean, $\hat{\mu}_n$, does not change each time it is computed like it would for IID nodes. While determinism has advantages, randomization has substantial advantages as well.
\begin{itemize}
	\item Randomization done correctly removes bias in the estimator, $\hat{\mu}_n$.
	\item Replications of random estimators can facilitate error estimates for the sample mean.
	\item Often the application of interest requires transforming the LD nodes to mimic other distributions.  This is the case in the Keister example of Section \ref{sec:intro}, as seen in \eqref{eq:sample_mean_Keister} where the nodes, $\{\bsx_i\}_{i=0}^{n-1}$, are transformed to mimic a Gaussian distribution.

	The LD sequences defined in the previous section start with $\bsx_0 = \bszero$, as can be seen in \eqref{eq:latticepts}, \eqref{eq:digital}, \eqref{eq:Halton}, and Figures \ref{fig:latticeconstruct}--\ref{fig:Halton}.  The node $\bszero$ becomes infinite under a transformation to mimic a Gaussian distribution, which can trigger runtime errors.  Randomization eliminates nodes on the boundary of the unit cube.
\end{itemize}
The key to good randomization is to preserve the LD quality of the node sequence.  In general, we recommend randomization, if available.

\subsection{Shifts} \label{sec:shifts}
\begin{figure}
	\centering
	\includegraphics[width = \textwidth]{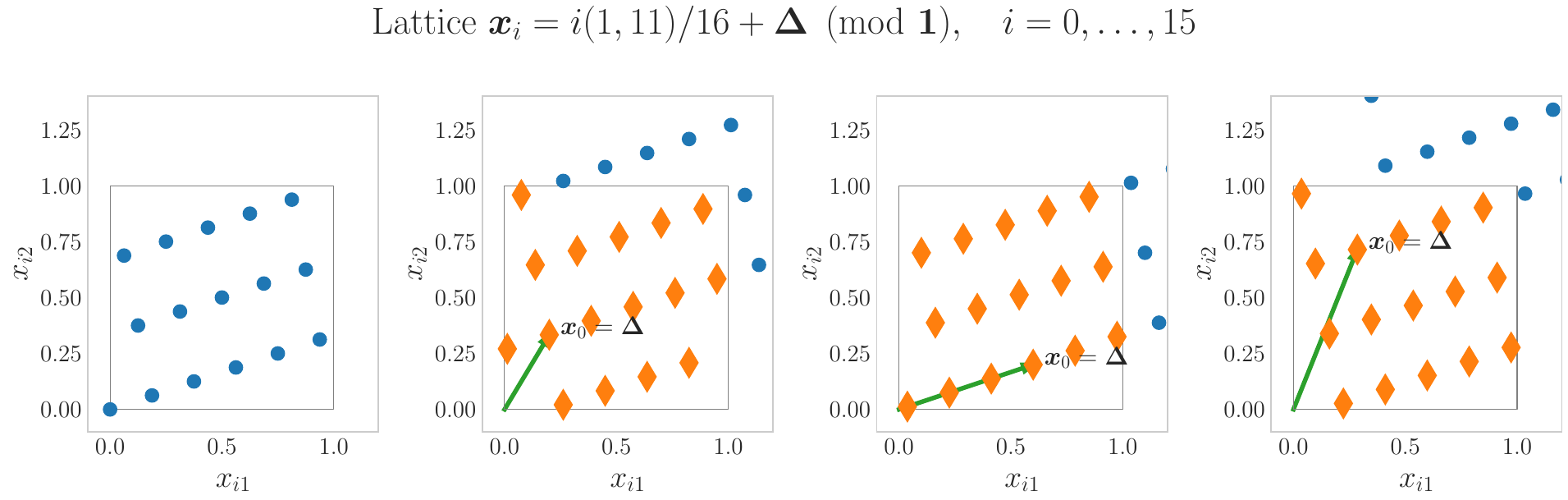}
	\caption{The original unshifted lattice (left) and three different shifts modulo $\bsone$ of the lattice.  Note that the structure of the lattice is maintained. \label{fig:shift_lat}}
\end{figure}

The simplest randomization is to shift the node sequence by $\bsDelta  \sim \calU[0,1)^d$.  For lattice  sequences the shift is applied modulo $\bsone$ so that shifted extensible lattice turns
\eqref{eq:extensiblelattice} into
\begin{equation} \label{eq:sh_extensiblelattice}
	\bsx^{\bsDelta\text{-lat}}_i := \bsx_i^{\text{lat}} + \bsDelta \pmod \bsone =\phi_b( i)\bsh + \bsDelta \pmod \bsone, \qquad i = 0, 1 , \ldots.
\end{equation}
Since the first $2^m$ elements of the lattice sequence form a group under modulo $\bsone$ addition, the corresponding shifted node set is a coset.  Figure \ref{fig:shift_lat} illustrates three shifts of the lattice plotted in Figure \ref{fig:latticeconstruct}.

\begin{figure}
	\centering
	\includegraphics[width = \textwidth]{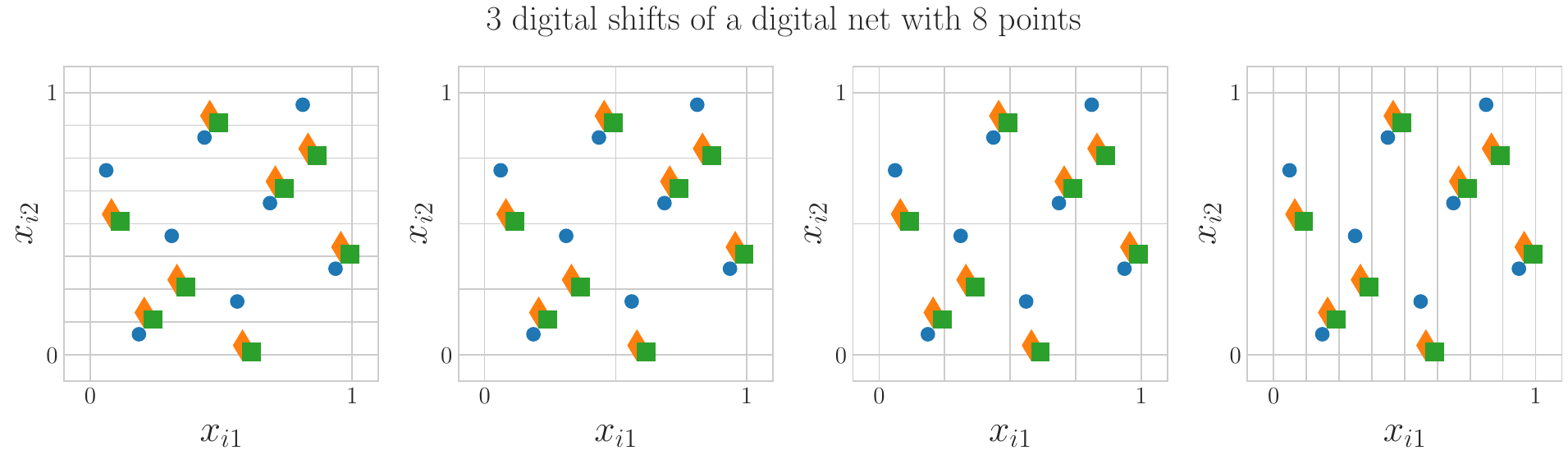}
	\caption{Three digital shifts of a digital $(0,3,2)$-net with nodes denoted by three different colors/shapes.  There is one node in each smaller box of the form \eqref{eq:element_box} for $\lVert \bsk \rVert = 3$, for each digitally shifted net. \label{fig:shift_net}}
\end{figure}

For digital sequences the shift should be applied using digitwise addition.  A digital shift of the digital sequence defined in  \eqref{eq:digital} then becomes
\begin{multline} \label{eq:sh_digital}
	\bsx_i^{\bsDelta\text{-dig}} := \bsx_i^{\text{dig}} \oplus \bsDelta = i_0 \bsx^{\text{dig}}_1 \oplus i_1 \bsx^{\text{dig}}_2 \oplus i_2 \bsx^{\text{dig}}_4 \oplus \cdots \oplus \bsDelta \\
    \text{for }
	i = i_0 + i_12 + i_2 2^2 + \cdots.
\end{multline}
Three digital shifts of the same net are given in Figure \ref{fig:shift_net}.  A digital shift does not alter the $t$-value of a $(t,m,d)$-net.

Let $\{\bsx_i^{\bsDelta}\}_{i=0}^{n-1}$ denote a uniform (digital or modulo $\bsone$) random shift of any deterministic original node set, $\{\bsx_i \}_{i=0}^{n-1}$.  That is, $\bsx_i^{\bsDelta} : = \bsx_i + \bsDelta \pmod \bsone$ or  $\bsx_i^{\bsDelta} : = \bsx_i \oplus \bsDelta$ and $\bsDelta \sim \calU[0,1)^d$. This then implies that each $\bsx_i^{\bsDelta} \sim \calU[0,1)^d$ and thus $\hat{\mu}_n$ is an unbiased estimator of $\mu$, which was mentioned as an advantage at the beginning of this section.

\subsection{Digital Scrambles} \label{sec:scrambles}
For digital nets one can randomize even further by a scrambling that preserves the $t$-value.  Owen \cite{Owe95} proposed a scrambling of nets called nested uniform scrambling.  A simpler version called linear matrix scrambling was proposed by Matou\v{s}ek \cite{Mat98} and is described here.  Starting from the formulation of digital sequences involving generator matrices in \eqref{eq:digitalB}, we multiply each generator matrix on the left by a lower triangular matrix with ones along the diagonals and elements below the diagonal that are randomly chosen to be $0$ or $1$ with equal probability:
\begin{subequations} \label{eq:scr_lattice}
\begin{gather} \label{eq:digitalscr}
	\begin{pmatrix} x^{\text{scr-net}}_{ij1} \\ x^{\text{scr-net}}_{ij2} \\ \vdots \end{pmatrix}
	= \mathsf{L_j} \mathsf{C}_j \begin{pmatrix} i_0 \\ i_1 \\ \vdots \end{pmatrix} + \begin{pmatrix} \Delta_{1j} \\ \Delta_{2j} \\ \vdots \end{pmatrix} \quad \pmod{2}, \qquad j = 1,\ldots, d, \ i = 0, 1, \ldots, \\
	\mathsf{L_j} :=
	\begin{pmatrix}
		1 & 0 & 0 & 0 & \cdots \\
		l_{21} & 1 & 0 & 0 & \cdots \\
		l_{31} & l_{32} & 1 & 0 & \cdots \\
		\vdots & \vdots & \vdots & \ddots & \ddots
	\end{pmatrix}, \qquad l_{\ell j} \overset {\text{IID}}{\sim} \calU\{0,1\}
\end{gather}
\end{subequations}
A random digital shift, $\bsDelta$ is also added, where $\Delta_{\ell j}$ denotes the $\ell^{\text{th}}$ digit of the $j^{\text{th}}$ dimension of the shift.  Digital scrambling of a $(t,m,d)$-net does not alter its $t$-value.

\begin{figure}
	\centering
	\includegraphics[width = \textwidth]{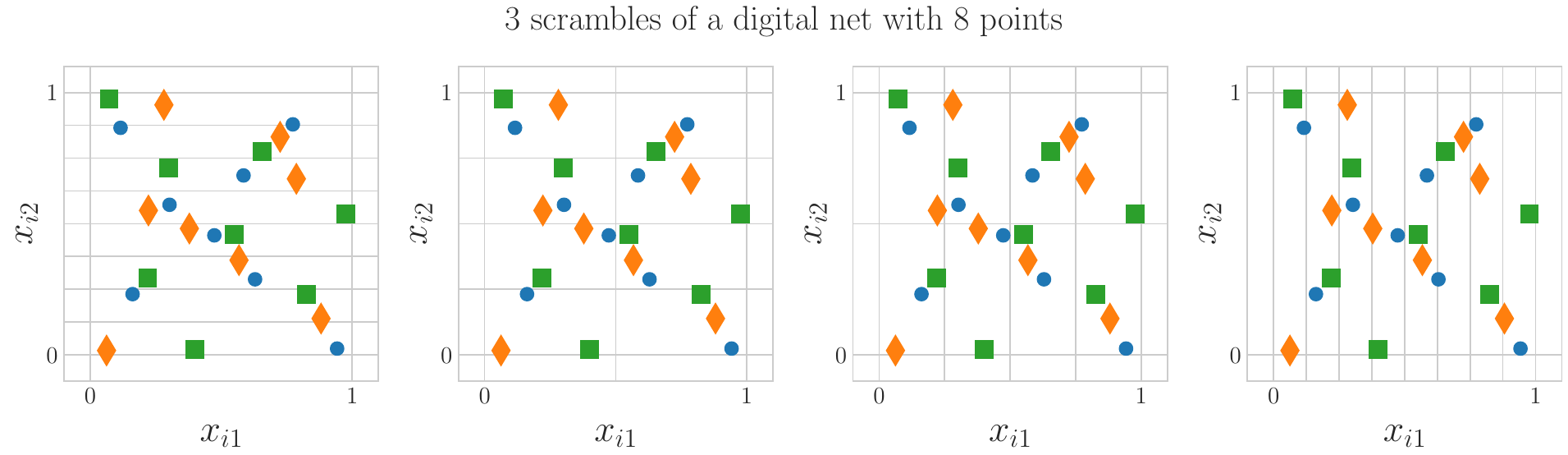}
	\caption{Three linear scrambles of a digital $(0,3,2)$-net with nodes denoted by three different colors/shapes.  There is one node in each smaller box of the form \eqref{eq:element_box} for $\lVert \bsk \rVert = 3$, for each digitally shifted net. \label{fig:scrambled_net}}
\end{figure}

\subsection{Randomizations for Halton}\label{sec:randHalton}


As Halton nodes can be viewed as digital nets with different bases for each dimension, the randomization techniques used for digital nets can be applied to Halton node sets. Random nested uniform scrambling and linear matrix scrambling have also been considered for Halton node sets \cite{owen2024gain}.

We detail one possible randomization of the Halton sequence which is equivalent to applying a random digital shift \eqref{eq:sh_digital} in the respective base of each dimension of the Halton sequence. We first define the generalized van der Corput sequence in base $b$ by obtaining a sequence $\Sigma = (\sigma_r)_{r\geq0}$ of permutations of $\{0,\ldots,b-1\}$ where the $i^{th}$ term of the one-dimensional sequence is defined as 
\begin{multline}\label{eq:gen_vdc}
    \phi_b^{\Sigma}(i_0 + i_1b + i_2 b^2 + \cdots ) := \sigma_0(i_0) b^{-1} + \sigma_1(i_1) b^{-2} + \sigma_2(i_2) b^{-3} + \cdots \in [0,1)
	\\
	 \text{where } i_0, \sigma_0(i_0), i_1, \sigma_1(i_1), \ldots \in \{0,\ldots b-1\}.
\end{multline}
The generalized Halton sequence is subsequently defined by choosing $d$ sequences of permutations $\Sigma_j = (\sigma_{j,r})_{r\geq0}$ of $\{0,\ldots,b_j-1\}$ for $j = 1, 2,\ldots,d$ and setting
$$\bsx_i^{\text{G-Hal}} = \left( \phi^{\Sigma_1}_{b_1}(i), \ldots, \phi^{\Sigma_d}_{b_d}(i)\right), \quad i= 0,1,\ldots,$$ A randomized qMC estimate $\hat{\mu}_n$ is then an unbiased estimator of $\mu$ if the permutations $\Sigma_j$ are sampled independently and uniformly, i.e $\mu$ if $\Sigma_j \overset{\text{IID}}{\sim} \calU\{0,\dots,b_j-1\}^d$ since this ensures $\bsx_i^{\text{G-Hal}} \sim \mathcal{U}[0,1)^d$. 

We note in passing that many works \cite{MasChiWar05,tuffin98,Atan04,vancools06} have sought to optimize the permutations $\Sigma_j$ deterministically to enhance uniformity. For a comprehensive survey of the Halton sequence (up to 2008), we refer readers to the excellent work by Faure and Lemieux \cite{faulem09}.

\subsection{Randomized Generators}

Instead of explicitly searching for good generating vectors or matrices, one can instead sample a few random generating vectors or generating matrices and then prune out bad choices. One recent method has explored selecting the median of $r$ randomized qMC estimates, rather than the mean \cite{PanOwe23a}.  Random generating vectors for lattices take $\bsh \overset{\text{IID}}{\sim} \calU\{1,\dots,n-1\}^d$ while random generating matrices for base $2$ digital nets---as described in \eqref{eq:digitalB}---set $x^{\text{dig}}_{ij\ell} \overset{\text{IID}}{\sim} \calU\{0,1\}$ for $i<\ell$ while enforcing $x_{iji}=1$ and $x_{ij\ell}=0$ for $i>\ell$. 

\section{Discrepancy} \label{sec:discrepancy}
So far, we have relied on eye tests and the Keister example to show how LD sequences are better than IID.  This section introduces the theory that shows why minimizing discrepancy leads to better approximations to $\mu$.

\subsection{Discrepancies Defined by Kernels} \label{sec:kerdisc}
Let $K: [0,1)^d \times [0,1)^d \to \mathbb{R}$ be a function satisfying the following two conditions:
\begin{subequations} \label{eq:Kcond}
	\begin{align}
		\label{eq:KcondSym}
		\text{Symmetry:} \quad & K(\bst,\bsx) = K(\bsx,\bst) \qquad \forall \bst,\bsx \in [0,1)^d \\
		\nonumber
		\text{Postitive definiteness:} \quad & \sum_{i,j = 0}^{n-1} c_i K(\bsx_i,\bsx_j) c_j > 0 \qquad  \forall n \in \mathbb{N}, \ \bsc \ne 0, \\
		& \qquad \qquad  \ \text{distinct } \bsx_0, \ldots, \bsx_{n-1} \in [0,1)^d. \label{eq:KcondPD}
	\end{align}
\end{subequations}
Then it is known \cite{Aro50} that this $K$ is the reproducing kernel for a Hilbert space, $\calH_K$ with associated inner product $\langle \cdot, \cdot \rangle_{\calH_K}$, such that
\begin{subequations}
	\begin{align}
	\text{Belonging:} \quad & K(\cdot,\bsx) \in \calH_K \qquad \forall \bsx \in [0,1)^d \\
	\text{Reproducing:} \quad & f(\bsx) = \langle K(\cdot,\bsx),f \rangle \qquad  \forall \bsx \in [0,1)^d, \ f \in \calH_K.
\end{align}
\end{subequations}

Having a reproducing kernel Hilbert space with a known kernel $K$ allows us to derive a rigorous error bound for $\mu - \hat{\mu}_n$.  First note that for any bounded, linear functional, $L$ on the Hilbert space $\calH_K$, there is some $\eta_L \in \calH_K$ such that
\begin{equation*}
L(f) = \langle \eta_L , f \rangle_{\calH_K} \qquad \forall f \in \calH_K.
\end{equation*}
This is guaranteed by the Riesz Representation Theorem, and $\eta_L$ is called the representer of $L$.  Using the reproducing property of $K$, one may derive an explicit formula for $\eta_L(\bsx)$, namely,
\begin{equation*}
\eta_L(\bsx) = \langle K(\cdot, \bsx), \eta_L \rangle_{\calH_K} = L\bigl( K(\cdot,\bsx) \bigr).
\end{equation*}
From this expression for $\eta_L(\bsx)$ one may then calculate the squared norm of the linear functional $L$ as the squared norm of its representer
\begin{equation} \label{eq:sqnormL}
	\lVert \eta_L \rVert_{\calH_K}^2 = \langle\eta_L, \eta_L \rangle_{\calH_K} = L(\eta_L) = L^{\cdot\cdot} \Bigl(L^{\cdot}\bigl( K(\cdot,\cdot \cdot) \bigr) \Bigr).
\end{equation}

If $\mu(\cdot):\calH_{\calK} \to \RR$ is a bounded linear functional on $\calH_K$, we may use the argument above to derive our error bound.  First we write the dependence of $\mu - \hat{\mu}_n$ on  $f$ explicitly as $\mu(f) - \hat{\mu}_n(f)$ and note that the error functional, $\mu(\cdot) - \hat{\mu}_n(\cdot)$, is  linear.  The Riesz Representation Theorem and the Cauchy-Schwarz inequality imply a tight error bound
\begin{align}
	\nonumber
	\lvert\mu(f) - \hat{\mu}_n(f) \rvert
	& = \lvert \langle \eta_{\mu(\cdot) - \hat{\mu}_n(\cdot)}, f \rangle_{\calH_K} \rvert
	 \le \lVert  \eta_{\mu(\cdot) - \hat{\mu}_n(\cdot)} \rVert_{\calH_K} \, \lVert f \rVert_{\calH_K}\\
	\lvert\mu(f) - \hat{\mu}_n(f) \rvert
	& \le \underbrace{\lVert  \eta_{\mu(\cdot) - \hat{\mu}_n(\cdot)} \rVert_{\calH_K}}_{=:\text{discrepancy}\bigl(\{\bsx\}_{i=0}^{n-1},K \bigr)}
	\, \underbrace{\inf_{c \in \RR}\lVert  f - c \rVert_{\calH_K}}_{=:\text{variation}(f,K)} \qquad \forall f \in \calH_K,  \label{eq:cuberrbd}
	\intertext{since $\hat{\mu}_n$ is exact for constants.  (It is assumed that $\calH_K$ contains constant functions so that $f \in \calH_K$ implies that $f - c \in \calH_K$.)  The squared norm of the error funtional can be expressed explicitly via \eqref{eq:sqnormL} as}
	\nonumber
	\text{discrepancy}^2\bigl(\{\bsx\}_{i=0}^{n-1},K \bigr) & = \lVert  \eta_{\mu(\cdot) - \hat{\mu}_n(\cdot)} \rVert_{\calH_K}^2 \\
	\nonumber
	& =
	\bigl(\mu(\cdot\cdot) - \hat{\mu}_n(\cdot\cdot)\bigr) \Bigl(\bigl(\mu(\cdot) - \hat{\mu}_n(\cdot)\bigr)\bigl( K(\cdot,\cdot \cdot) \bigr) \Bigr) \\
	\nonumber
	& = \int_{[0,1)^d \times [0,1)^d} K(\bst,\bsx) \, \mathrm{d} \bst \, \mathrm{d} \bsx  -\frac 2n  \sum_{i=0}^{n-1} \int_{[0,1)^d} K(\bst,\bsx_i) \, \mathrm{d} \bst \\
	& \qquad \qquad + \frac{1}{n^2} \sum_{i,j=0}^{n-1}  K(\bsx_i,\bsx_j). \label{eq:sqdisc}
\end{align}

The error bound \eqref{eq:cuberrbd} has two factors:
\begin{itemize}
	\item The discrepancy, which is the norm of the error functional, depends only on the node set, $\{\bsx\}_{i=0}^{n-1}$, and \emph{not} on $f$.  The discrepancy measures the deficiency of that node set.
	\item The variation depends only on $f$, the function defining our random variable $Y$ whose expectation we wish to compute,  and \emph{not} on the node set, $\{\bsx\}_{i=0}^{n-1}$.  The variation is a measure of the roughness of $f$.
\end{itemize}
In general, designers of qMC methods want to construct sets or sequences of nodes with  discrepancy as small as possible, either by increasing $n$, if the computational budget allows, or by better placement of the nodes.  The constructions described in Section \ref{sec:construct} are ways of ensuring better placement of the nodes, provided that the generators are chosen well.  Practitioners of qMC want to formulate $\mu$ in a way to make the variation as small as possible.  This is discussed in Section \ref{sec:reformulate}.

Although error bound \eqref{eq:cuberrbd} is elegant, it leaves several matters unresolved:
\begin{itemize}
	\item The reproducing kernel $K$ defines both the discrepancy and the variation.  This means that different choices of $K$ will lead to different error bounds for the same $\hat{\mu}_n$, even though the error is unchanged.
	\item The Hilbert space $\calH_K$ contains one's assumptions about $f$, such as smoothness or periodicity. Knowing $K$, however, does not automatically lead to an explicit formula for the inner product of the associated Hilbert space, $\langle \cdot, \cdot \rangle_{\calH_K}$.  Conversely, having an explicit formula for  $\langle \cdot, \cdot \rangle_{\calH_K}$ does not automatically lead to an explicit formula for the reproducing kernel, $K$.  In both cases, educated guesswork is involved.
	\item Choosing a more restrictive $\calH_K$ may lead to a sharper error bound, however, it is often infeasible in practice to check whether the $f$ lies in $\calH_K$.  This is what we mean in the introduction by referring to $f$ as a ``black box''.
	\item   Moreover, even if one is confident that $f \in \calH_K$ ($f$ is a ``gray box''), it is usually impractical to compute variation$(f,K)$, as we shall see in the discrepancy example below.  Error bound \eqref{eq:cuberrbd} cannot be used as a stopping criterion to determine the $n$ needed to satisfy an error tolerance.
\end{itemize}
In spite of these drawbacks, \eqref{eq:cuberrbd} is quite useful:
\begin{itemize}
	\item There are families of commonly used reproducing kernels \cite{BerT-A04,Hic97a,Hic98b,Hic99b}. Once $K$ has been fixed, the discrepancy may be used to compare different LD sequences.  The discrepancy's rate of decay---which sometimes can be computed theoretically---indicates how quickly $\hat{\mu}_n$ approaches $\mu$.
	\item When an explicit formula for the inner product of the associated Hilbert space, $\langle \cdot, \cdot \rangle_{\calH_K}$, can be deduced, it suggests how to formulate $\mu$ to reduce variation$(f,K)$.
\end{itemize}

\subsection{The Centered Discrepancy}
An example of a reproducing kernel defined on $[0,1)^d \times [0,1)^d$ is
\begin{equation} \label{eq:ctrkernel}
	K(\bst,\bsx) = \prod_{\ell = 1}^d \left[ 1 + \frac 12 \bigl ( \lvert t_\ell - 1/2 \rvert + \lvert x_\ell - 1/2 \rvert - \lvert t_\ell - x_\ell \rvert \bigr ) \right] .
\end{equation}
One can verify that this $K$ satisfies the symmetry condition \eqref{eq:KcondSym}.  It also is symmetric about the middle of the unit cube, i.e.,
\begin{equation*}
	K\bigl((t_1,\ldots, t_{\ell-1}, 1 - t_\ell,  t_{\ell+1}, \ldots, t_d),
	(x_1,\ldots, x_{\ell-1}, 1 - x_\ell,  x_{\ell+1}, \ldots, x_d)
	\bigr) =
	K(\bst,\bsx),
\end{equation*}
which is a desirable property for problems where there is no preferred corner of the unit cube.

The inner product for the Hilbert space $\calH_K$ with reproducing kernel $K$ is
\begin{align}
	\langle f, g \rangle_{\calH_K} &= f(1/2, \ldots, 1/2) g(1/2, \ldots, 1/2) \nonumber
	\\
	&
	+ \int_{[0,1)} [D^{\{1\}}f](x_1) \, [D^{\{1\}}g](x_1)\, \mathrm{d} x_1
	+ \int_{[0,1)} [D^{\{2\}}f](x_2) \, [D^{\{2\}}g](x_2)\, \mathrm{d} x_2 \nonumber \\
	& \qquad \qquad + \cdots \nonumber\\
	&
	+ \int_{[0,1)^2} [D^{{\{1,2\}}}f](x_1,x_2) \, [D^{\{1,2\}}g](x_1,x_2)\, \mathrm{d} x_1 \mathrm{d} x_2 + \cdots \nonumber \\
	&  \qquad \qquad + \cdots + \int_{[0,1)^d} [D^{\{1,\ldots, d\}}f](\bsx) \, [D^{\{1,\ldots, d\}}g](\bsx) \, \mathrm{d} \bsx \nonumber \\
	& = \sum_{\bsell \subseteq \{1, \ldots, d\}} \langle D^{\bsell}f, D^{\bsell}g \rangle_2
\end{align}
where $\langle \cdot,\cdot \rangle_2$ denotes the $L^2$ inner product, the operator $D^{\bsell}$ is defined as
\begin{multline*}
	[D^{\{\ell_1, \ell_2, \ldots, \ell_s\} }f](x_1, x_2, \ldots, x_s) \\
	: = \frac{\partial^s f(1/2,\ldots, 1/2, x_{\ell_1}, 1/2, \ldots, 1/2, x_{\ell_2}, 1/2, \ldots, 1/2, x_{\ell_s}, 1/2, \ldots, 1/2)}{\partial x_{\ell_1} \partial x_{\ell_2} \cdots \partial x_{\ell_s}},
\end{multline*}
and $D^{\emptyset}f$ is understood to be the constant $f(1/2, \ldots, 1/2)$ \cite{Hic97a}.  The variation, which corresponds to the non-constant part of the function, is
\begin{equation*}
		\text{variation}(f) := \sqrt{\sum_{\emptyset \ne \bsell \subseteq \{1, \ldots, d\}} \lVert D^{\bsell}f\rVert_2^2}.
\end{equation*}

We show that $K$ defined in \eqref{eq:ctrkernel} is the reproducing kernel for the $\calH_K$ with the above inner product for $d=1$ and assuming $0 \le x \le 1/2$:
\begin{align*}
	\langle K(\cdot,x), f \rangle_{\calH_K}
	& =
	K(1/2,x) f(1/2) + \int_0^1 \frac{\partial K(t,x)}{\partial t} f'(t) \, \mathrm{d} t \\
	& = 1 \times f(1/2) +  \int_0^x \frac 12 [-1  - (-1) ] f'(t) \, \mathrm{d} t +
	\int_x^{1/2} \frac 12 [-1  - 1 ] f'(t) \, \mathrm{d} t \\
	& \qquad \qquad +
	\int_{1/2}^{1} \frac 12 [1  - 1 ] f'(t) \, \mathrm{d} t \\
	& = f(1/2) + 0 - \int_x^{1/2} f'(t) \, \mathrm{d} t + 0 = f(1/2) - [f(1/2) - f(x)] \\
	& =  f(x).
\end{align*}
The extensions to general $1/2 < x < 1$ and to $d > 1$ are straightforward \cite{Hic97a}.

\begin{figure}
	\centering
	\includegraphics[width = 0.48\textwidth]{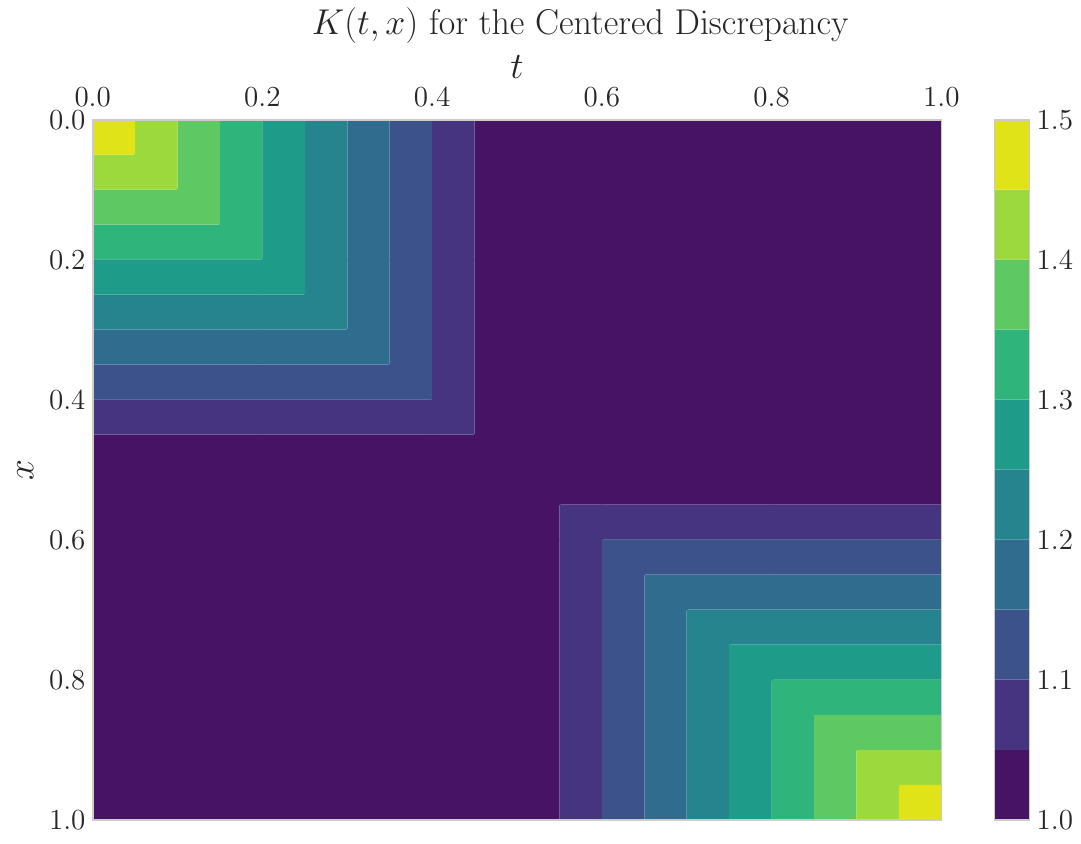} \quad 
    \includegraphics[width = 0.48\textwidth]{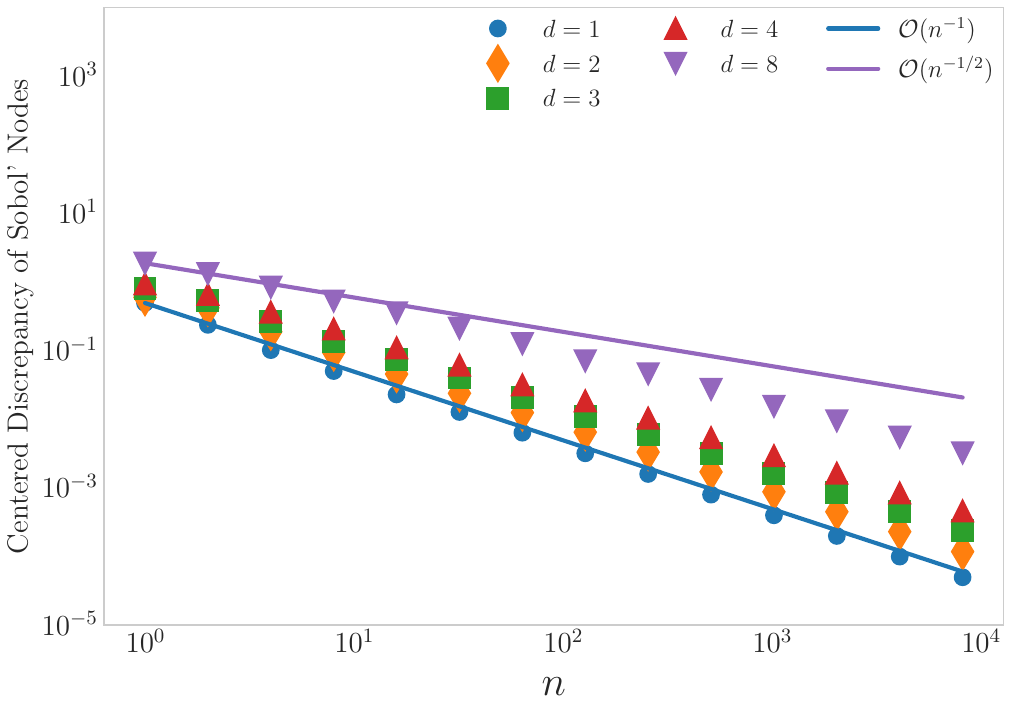}
	\caption{Left:  A plot of the reproducing kernel for the centered discrepancy defined in \eqref{eq:ctrkernel} for $d=1$.  Note that the vertical axis is inverted so that the plot is a picture of the ``matrix'' $\bigl( K(t,x) \bigr)_{t,x \in [0,1)]}$.  Any submatrix is positive definite by property \eqref{eq:KcondPD}. \\
    Right: The centered discrepancy for several modest values of $d$.  They show an asymptotic decay close to $\calO(n^{-1})$, but this rate is achieved only at higher $n$ for larger values of $d$ \label{fig:ctrKer}}
\end{figure}

Figure \ref{fig:ctrKer} (left) is a color map of the reproducing kernel, $K$, defined in \eqref{eq:ctrkernel} for $d=1$.  As shown in the figure, the matrix $\mathsf{K} := \bigl( K(x_i,x_j) \bigr)_{i,j=1}^n$ that arises in condition \eqref{eq:KcondPD} is formed by taking $n$ rows and corresponding columns of this plot.  Any such matrix $\mathsf{K}$ is strictly positive definite, as is suggested by the higher values near the diagonal.

The squared centered discrepancy is given by substituting the formula for $K$ into \eqref{eq:sqdisc}, which implies
\begin{align} \label{eq:ctrdisc}
	\nonumber
	\MoveEqLeft\text{ctr-discrepancy}^2\bigl(\{ \bsx_i\}_{i=0}^{n-1} \bigr) \\
	\nonumber
	&= \left ( \frac{13}{12} \right )^d
	- \frac 2n \sum_{i=0}^{n-1} \prod_{\ell=1}^d \left [ 1 + \frac 12 \bigl ( \lvert x_{i\ell} - 1/2 \rvert - \lvert x_{i\ell} - 1/2 \rvert^d \bigr )\right ] \\
	& \qquad \qquad \frac{1}{n^2} \sum_{i,j=0}^{n-1} \prod_{\ell = 1}^d \left[ 1 + \frac 12 \bigl ( \lvert x_{i\ell} - 1/2 \rvert + \lvert x_{j\ell} - 1/2 \rvert - \lvert x_{i\ell} - x_{j\ell} \rvert \bigr ) \right].
\end{align}
As is true for most cases, the computational cost of computing the discrepancy is $\mathcal{O}(dn^2)$ as $d$ and/or $n$ tend to infinity.

Figure  \ref{fig:ctrKer} (right)  displays the discrepancy for a Sobol' LD digital net for several modest values of $d$ and for $n =  1, 2, 4, 8, 16, \ldots$.  The asymptotic decay rate of nearly $\calO(n^{-1})$ is observed for small $d$ and can be anticipated for a bit larger $d$.  Thus, the decay of the discrepancy mirrors the convergence rate for LD sequences for the Keister example in Section \ref{sec:intro} and Figure \ref{fig:keister-err}.

The centered discrepancy also has a geometric interpretation, as is described in \cite{Hic97a}. Other choices of kernels and their geometric interpretations are explained there as well.  Moreover, it is possible to define discrepancies for Banach spaces of integrands, however, the formulas for these discrepancies are not as amenable to calculation.

\subsection{Coordinate weights} \label{sec:coordwts}
\begin{figure}
    \centering
     \includegraphics[width=0.45\textwidth]{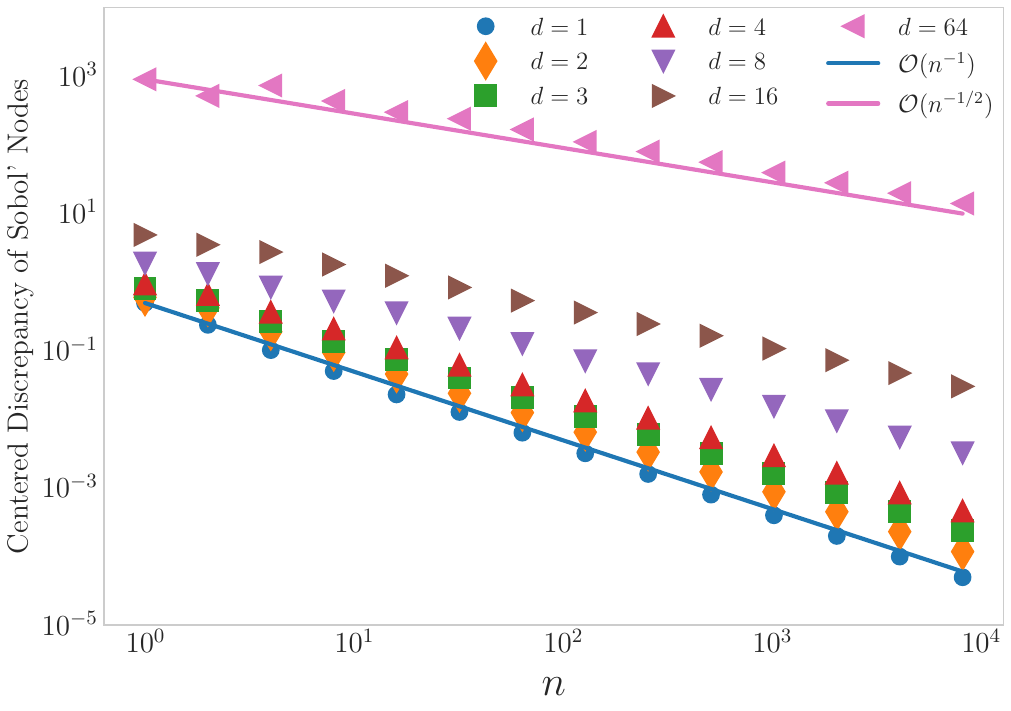}
    \includegraphics[width=0.45\textwidth]{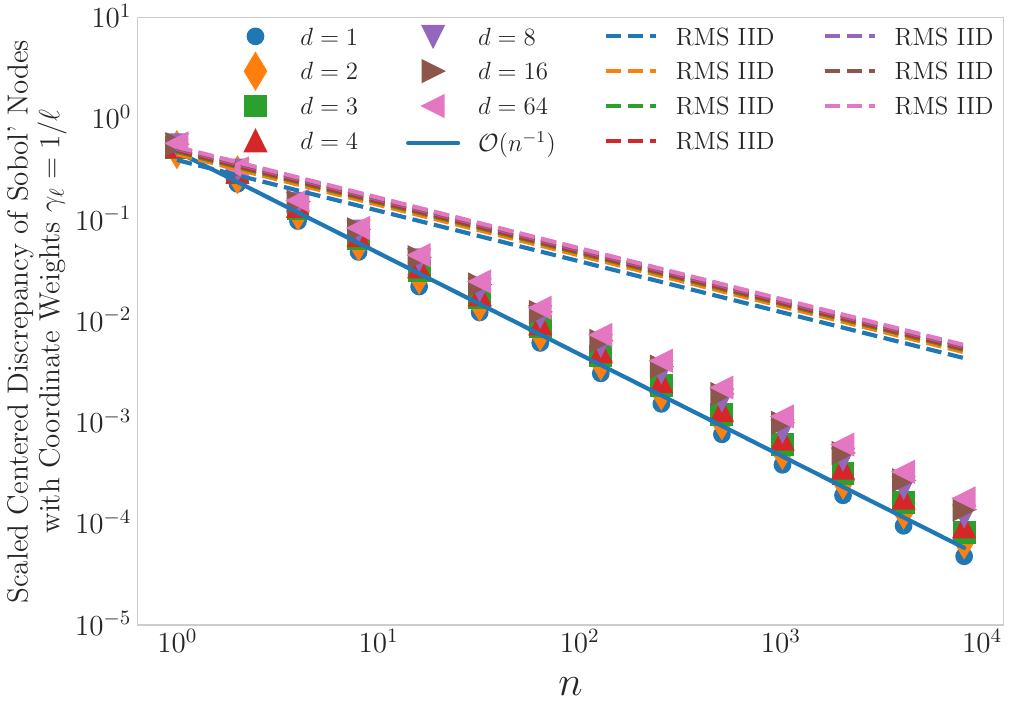}
    \caption{Left: The centered discrepancy for a Sobol' LD sequence for a range of $d$.  For larger $d$ only a decay rate of $\calO(n^{-1/2})$, is achieved for the range of $n$ shown.  Moreover, the discrepancy for larger $d$ is huge.
    Right:
    The scaled centered discrepancy with coordinate weights for a Sobol' LD sequence for various $d$ and the root mean squared scaled centered discrepancy with coordinate weights for IID nodes.  LD nodes have a much smaller discrepancy than IID nodes.}
    \label{fig:scaledweightedcentered}
\end{figure}

When the dimension, $d$, is increased, the discrepancy increases and the rate of decay with increasing sample size deteriorates, as can be seen in Figure \ref{fig:scaledweightedcentered} (left). This raises the question of whether the LD sequences lose their effectiveness for higher dimensions or whether the theory of the previous subsection does not coincide with practice.

Note that the discrepancy of the null set, which is equivalent to the norm of the linear functional $\mu(\cdot)$, is
\begin{equation}
		\text{discrepancy}\bigl(\emptyset,K \bigr)  = \sqrt{\int_{[0,1)^d \times [0,1)^d} K(\bst,\bsx) \, \mathrm{d} \bst \, \mathrm{d} \bsx},
\end{equation}
which corresponds to $(13/12)^{d/2}$ for the centered discrepancy.  This means that integration gets exponentially harder as $d$ increases.  One may scale the discrepancy of a node set by dividing it by discrepancy$(\emptyset,K)$, but this does not recover the desired nearly $\mathcal{O}(n^{-1})$ decay rate.

A good approach is to introduce coordinate weights, $\gamma_1, \ldots, \gamma_d$ into the reproducing kernel.  These were introduced in \cite{Hic95} but begin to be fully exploited in \cite{SloWoz98} and the articles that followed.  For the kernel defining the centered discrepancy, this means
\begin{align} \label{eq:wtctrkernel}
	K_{\bsgamma}(\bst,\bsx) &= \prod_{\ell = 1}^d \biggl [ 1 + \frac {\gamma_\ell^2}2 \bigl ( \lvert t_\ell - 1/2 \rvert + \lvert x_\ell - 1/2 \rvert - \lvert t_\ell - x_\ell \rvert \bigr ) \biggr] , \\
	\langle f, g \rangle_{\calH_{K_{\bsgamma}}}  &= \sum_{\bsell \subseteq \{1, \ldots, d\}} \biggl(\prod_k \gamma_{\ell_k}^2 \biggr) \langle D^{\bsell}f, D^{\bsell}g \rangle_2, \\
    \label{eq:wtctrvariation}
	\text{variation}(f,K_{\bsgamma})  &= \sqrt{\sum_{\bsell \subseteq \{1, \ldots, d\}} \biggl( \prod_k \gamma_{\ell_k}^{-2} \biggr) \lVert D^{\bsell}f \rVert_2^2}, \\
	\label{eq:wt-ctrdisc}
	\nonumber
	\MoveEqLeft[4]\text{wt-ctr-discrepancy}^2\bigl(\{ \bsx_i\}_{i=0}^{n-1} , \bsgamma\bigr) \\
	\nonumber
	&= \prod_{\ell=1}^d \left ( 1 + \frac{\gamma_{\ell}^2}{12} \right )^d
	- \frac 2n \sum_{i=0}^{n-1} \prod_{\ell=1}^d \left [ 1 + \frac {\gamma_{\ell}^2}2 \bigl ( \lvert x_{i\ell} - 1/2 \rvert - \lvert x_{i\ell} - 1/2 \rvert^d \bigr )\right ] \\
	& \qquad \qquad \frac{1}{n^2} \sum_{i,j=0}^{n-1} \prod_{\ell = 1}^d \left[ 1 + \frac {\gamma_{\ell}^2}2 \bigl ( \lvert x_{i\ell} - 1/2 \rvert + \lvert x_{j\ell} - 1/2 \rvert - \lvert x_{i\ell} - x_{j\ell} \rvert \bigr ) \right].
\end{align}

Figure \ref{fig:scaledweightedcentered} (right) plots the scaled centered discrepancy with coordinate weights $\bsgamma = (1/\ell)_{\ell =1}^d$, i.e., 
\begin{equation*}
    \text{wt-ctr-discrepancy}\bigl(\{ \bsx_i\}_{i=0}^{n-1},\bsgamma \bigr)/\text{wt-ctr-discrepancy}(\emptyset,\bsgamma ),
\end{equation*}
against sample size for a variety of $d$.  Now, even for larger $d$ there is a reasonable decay that approximates $\mathcal{O}(n^{-1})$.

For IID samples the root mean squared discrepancy for a general reproducing kernel, $K$, is
\begin{multline} \label{eq:discIID}
    \sqrt{\mathbb{E}[\text{discrepancy}^2\bigl(\{\bsx\}_{i=0}^{n-1},K \bigr)] } \\
    = \frac 1{\sqrt{n}}  \sqrt{  \int_{[0,1)^d} K(\bsx,\bsx) \, \mathrm{d} \bsx -
	 \int_{[0,1)^d \times [0,1)^d} K(\bst,\bsx) \, \mathrm{d} \bst \, \mathrm{d} \bsx}.
\end{multline}
This is also plotted in Figure \ref{fig:scaledweightedcentered} (right). which highlights the superiority of LD nodes.

The meaning of the coordinate weights is that those coordinates with a small weight should not contribute much to the variation.  In \eqref{eq:wtctrvariation}, if the $\gamma_{\ell_1}, \ldots, \gamma_{\ell_s}$ are small, then $\prod_k \gamma_{\ell_k}^{-2}$ will be large.  For $f$ to have modest variation requires $\lVert D^{\bsell}f \rVert_2$ to be small.  Even though a function may have a large nominal dimension, $d$, its effective dimension, i.e., the number of variables that substantially contribute to the variation, must be relatively few, so that the variation under coordinate weights is modest.  Under these conditions, the corresponding discrepancy may retain the nearly $\mathcal{O}(n^{-1})$ decay.

\subsection{Discrepancies for Lattices and Digital Nets} \label{sec:latticenetdis}
Many useful constructions of LD sequences are constructed numerically by minimizing the discrepancy, however, as noted earlier, computing the discrepancy requires $\mathcal{O}(dn^2)$ operations, which can be expensive.  However, this computational cost can be reduced for lattices and digital nets when using friendly kernels.

The mean squared discrepancy of a randomly shifted node set, $\{\bsx_i^{\bsDelta} =  \bsx_i + \bsDelta \bmod \bsone \}_{i=0}^{n-1}$ for an arbitrary kernel takes the form
\begin{align}
\nonumber
\MoveEqLeft \mathbb{E}\bigl[\text{discrepancy}^2 \bigl (\{ \bsx^{\bsDelta}_{i} \}_{i=0}^{n-1},K \bigr) \bigr] \\
\nonumber
& = \mathbb{E}\bigl[\text{discrepancy}^2 \bigl (\{ \bsx_{i} + \bsDelta \bmod \bsone \}_{i=0}^{n-1},K \bigr) \bigr] \\
\nonumber
& = \int_{[0,1)^d \times [0,1)^d} K(\bst,\bsx) \, \mathrm{d} \bst \, \mathrm{d} \bsx  -\frac 2n  \sum_{i=0}^{n-1}  \int_{[0,1)^d} \mathbb{E} \bigl [K(\bst,\bsx_i + \bsDelta \bmod \bsone) \bigr]  \, \mathrm{d} \bst \\
\nonumber
	& \qquad \qquad + \frac{1}{n^2} \sum_{i,j=0}^{n-1}  \mathbb{E} \bigl [ K(\bsx_i + \bsDelta \bmod \bsone,\bsx_j + \bsDelta \bmod \bsone) \bigr] \\
\nonumber
    & = \frac{1}{n^2} \sum_{i,j=0}^{n-1}  \mathbb{E} \bigl [ K(\bsx_i - \bsx_j  + \bsDelta \bmod \bsone, \bsDelta \bmod \bsone) \bigr] - \int_{[0,1)^d \times [0,1)^d} K(\bst,\bsx) \, \mathrm{d} \bst \, \mathrm{d} \bsx ,\\
    \intertext{where we have applied a variable transformation $\bsDelta \to \bsDelta - \bsx_j$,}
    & = \frac{1}{n^2} \sum_{i,j=0}^{n-1}  \widetilde{K}(\bsx_i - \bsx_j \bmod \bsone ) - \int_{[0,1)^d} \widetilde{K}(\bsx) \, \mathrm{d} \bsx, \\
    \intertext{and where the filtered version of the reproducing kernel is defined as}
    \widetilde{K}(\bsx) &: = \int_{[0,1)^d} K(\bsx + \bst \bmod \bsone, \bst) \, \mathrm{d} \bst.
\end{align}

Note that 
\begin{equation} \label{eq:shinvK}
    K^{\text{shift}} (\bst, \bsx) := \widetilde{K}(\bst - \bsx \bmod \bsone)
\end{equation}
is a reproducing kernel.  We call it a shift-invariant kernel because $K^{\text{shift}} (\bst + \bsDelta \bmod \bsone, \bsx + \bsDelta \bmod1) = \widetilde{K}(\bst - \bsx \bmod \bsone) =  K^{\text{shift}} (\bst, \bsx) $ for all $\bsDelta \in [0,1)^d$.

For arbitrary node sets, this mean squared discrepancy is of similar computational cost as the original, but if the unshifted node set $\{\bsx_i\}_{i=0}^{n-1}$ is a lattice, $\{\bsx^{\text{lat}}_i\}_{i=0}^{n-1}$, as defined in \eqref{eq:latticepts}, then for any $j$, the set $\{\bsx^{\text{lat}}_i - \bsx^{\text{lat}}_j \bmod \bsone \}_{i=0}^{n-1}$ has the same values for all $j = 0, \ldots, n-1$.  Thus the mean squared discrepancy for a randomly shifted lattice can be computed in only $\mathcal{O}(dn)$ operations, assuming that $K(\bsx)$ requires $d$ operations to compute:
\begin{multline}
\mathbb{E}\bigl[\text{discrepancy}^2 \bigl (\{ \bsx^{\bsDelta-\text{lat}}_{i} \}_{i=0}^{n-1}, K \bigr) \bigr] = \frac{1}{n} \sum_{i=0}^{n-1}  \widetilde{K}(\bsx_i) - \int_{[0,1)^d} \widetilde{K}(\bsx) \, \mathrm{d} \bsx \\
= \text{discrepancy}^2 \bigl (\{ \bsx^{\text{lat}}_{i} \}_{i=0}^{n-1}, K^{\text{shift}} \bigr).
\end{multline}
For the kernel defining the centered discrepancy with coordinate weights,
\begin{equation}
    \widetilde{K}(\bsx) = \int_{[0,1)^d} \prod_{\ell = 1}^d \Bigl [ 1 + \gamma_\ell^2 \{ 1/4 - x_\ell (1-x_\ell ) \} \Bigr], \quad \int_{[0,1)^d} \widetilde{K}(\bsx) \, \mathrm{d} \bsx = \prod_{\ell = 1}^d \biggl [ 1 + \frac{\gamma_\ell^2}{12} \biggr ].
\end{equation}

Under modest smoothness assumptions, the discrepancy decays nearly as fast as $\mathcal{O}(n^{-1})$ for typical LD sequences.  However, when the reproducing kernel, $K$, is constructed to ensure that functions in $\calH_K$ have mixed partial derivatives of up to order $r$ in each coordinate and are periodic in all lower order derivatives, then the discrepancy of lattices can achieve a nearly $\mathcal{O}(n^{-r})$ decay rate.

For digital sequences, there is an analogous development of filtering the original kernel by digitally shifts and/or scrambling.  This provides a root mean square discrepancy that is again requires only $\mathcal{O}(dn)$ operations to compute, not $\mathcal{O}(dn^2)$ operations.

\subsection{Randomized Error and Bayesian Error} \label{sec:randBayes}
The error analysis leading to error bound \eqref{eq:cuberrbd} is a worst-case error analysis, meaning that the error in approximating $\mu(f)$ by $\hat{\mu}_n(f)$ is definitely no larger than the discrepancy times the variation, provided that $f$ lies inside the Hilbert space with reproducing kernel $K$.  We briefly mention other error measures here.  More detail is given in \cite{HicEtal17a}.

For sample means where the nodes, $\{\bsx_i\}_{i=0}^{n-1}$ are randomized, one can define a measure of their quality based on the root mean squared error of $\hat{\mu}_n(f)$:
\begin{align}
    \text{RMSE}(\hat{\mu}_n(f)) & : = \sqrt{\mathbb{E}[(\mu(f) - \hat{\mu}_n(f))^2]} \\
    &\le \underbrace{ \sup_{\text{variation}(f,K) \le 1} \sqrt{\mathbb{E}[(\mu(f) - \hat{\mu}_n(f))^2]}}_{\text{RERR}\bigl(\{\bsx_i\}_{i=0}^{n-1},K\bigr)} \, \text{variation}(f,K).
\end{align}
This quality measure $\text{RERR}\bigl(\{\bsx_i\}_{i=0}^{n-1},K\bigr)$ is typically more optimistic (smaller), than the root mean squared of  $\text{discrepancy}\bigl(\{\bsx_i\}_{i=0}^{n-1},K\bigr)$ because the latter finds a worst case function for each random instance of nodes, while the former is finds just one worst function.

Unfortunately, $\text{RERR}\bigl(\{\bsx_i\}_{i=0}^{n-1},K\bigr)$ is typically not as easy to compute as the discrepancy.  If $\hat{\mu}_n$ is unbiased, then $\sqrt{\mathbb{E}[(\mu(f) - \hat{\mu}_n(f))^2]} = \text{std}(\hat{\mu}_n(f))$.  For IID nodes we know that this is $\text{std}(f(\bsX))/\sqrt{n}$, which has the same order of decay as root mean squared discrepancy of an IID node set in \eqref{eq:discIID}, but a smaller constant.

For randomly shifted lattices $\text{RERR}\bigl(\{\bsx_i\}_{i=0}^{n-1},K\bigr)$ tends to have the same rate of decay as does the root mean squared discrepancy.  However, for randomly scrambled and shifted digital nets, there can be an improvement in rate of decay assuming a bit more smoothness \cite{Owe97,HeiHicYue02a,HicYue00}.

The symmetric, positive definite reproducing kernel, $K$, which appears in the worst-case error analysis in Section \ref{sec:kerdisc} may instead take the role of the covariance kernel of a Gaussian process with zero mean.  Average-case error analysis \cite{Rit00a} and probabilistic numerics \cite{BriEtal18a} pursue this approach.

\section{Stopping Criteria} \label{sec:stop}


When approximating the true $\mu$ with $\hat{\mu}_n$, stopping criterion determine how large $n$ should be so that the error $\lvert \mu - \hat{\mu}_n \rvert$ is below an error tolerance $\varepsilon$ either certainly or with sufficiently low uncertainty $\alpha$. The error tolerance and probability threshold are problem dependent. For example, in option pricing one may set the error tolerance to a penny, $\varepsilon=0.01$, and wish to find error bounds which hold with 95\% confidence, $\alpha=0.05$. The stopping criteria presented here are implemented in QMCPy \cite{QMCPy2020a}.

For simple Monte Carlo with IID nodes, the Central Limit Theorem (CLT) may be used to get an approximate confidence interval $[\mu^-,\mu^+]$ where $\mu^\pm = \hat{\mu} \pm Z^{1-\alpha/2} \sigma/\sqrt{n}$, $\sigma = \mathrm{Std}[f(\bsX)]$ and $Z^{1-\alpha/2}$ is the $1-\alpha/2$ quantile of the standard Gaussian distribution. If 
\begin{equation}
    n > (2 Z^{1-\alpha/2} \sigma/\varepsilon)^2
    \label{eq:sc_clt_iid_n}
\end{equation} then $\mu \in [(]\mu^-,\mu^+]$ with probability $1-\alpha$. A simple IID stopping criterion uses an initial sample of $n_0$ nodes to approximate $\sigma^2 \approx S^2 := (n_0-1)^{-1} \sum_{i=0}^{n_0-1}(f(\bsx_i)-\hat{\mu})^2$, and then plugs this into \eqref{eq:sc_clt_iid_n} to determine $n$. (Sometimes the sample standard deviation, $S$, is multiplied by an inflation factor to account for the fact that it is an unbiased estimate and not an upper bound.) 

As the CLT holds only as $n$ goes to infinity, this is an approximate stopping criterion. For integrands $f$ with bounded kurtosis, a guaranteed version of this two-stage-sampling idea has been developed to account for the finite sample size using  Berry-Esseen inequalities\cite{HicEtal14a}. 

A straightforward qMC stopping criterion may also be derived from the CLT using several randomizations of a LD node set. Let $\bsx_{ir}$ be the $i^\mathrm{th}$ node in the $r^\mathrm{th}$ randomization of a LD node set with $0 \leq i < n$ and $1 \leq r \leq R$. For example, we may use $R$ independent shifts of a lattice node set or independent digital shifts of a digital net. Crucially, the randomized node sets presented in Section \ref{sec:random} are also LD. Let $\hat{\mu}^{(r)} = n^{-1} \sum_{i=0}^{n-1} f(\bsx_{ir})$ be sample means for each randomization and set $\hat{\mu}_{nR} = R^{-1} \sum_{r=1}^R \hat{\mu}^{(r)}$ to be the aggregate approximation. In the spirit of CLT, a standard approximate $1-\alpha$ confidence interval $(\mu^-,\mu^+)$ has $\mu^\pm = \hat{\mu}_{nR} \pm t_{R-1}^{1-\alpha/2} S/\sqrt{R}$ where $S^2 = (R-1)^{-1} \sum_{r=1}^R (\hat{\mu}^{(r)} - \hat{\mu}_{nR})^2$ and $t_{R-1}^{1-\alpha/2}$ is the $1-\alpha/2$ quantile of the Student $t$ distribution with $R-1$ degrees of freedom. Here we have used the Student $t$ distribution since $R$ is typically small, for example $R=15$. To use this confidence interval for a stopping criterion, one may iteratively double the sample size $n$ until the interval width is sufficiently low. In \cite{l2023confidence}, comprehensive numerical experiments show that this straightforward qMC stopping criterion achieve better coverage than more involved bootstrap $t$ methods across a variety of test functions. 

One of the drawbacks of this qMC-CLT stopping criterion is that, at least intuitively, independent randomizations of LD nodes are less efficient than simply using more LD nodes from a single sequence. For instance, we would expect an approximation with $n=2^8$ and $R=2^3$ to perform worse than simply taking $n=2^{11}$ nodes from a single randomized LD sequence ($R=1$), but we need $R$ large enough to believe the CLT stopping criterion. 

One idea for constructing confidence intervals without replications is to approximate the decay rate of Fourier coefficients of $f$. Bounds of absolute certainty have be derived for cones of functions with reasonable decay of their Fourier coefficients \cite{HicEtal17a}. As in the qMC-CLT stopping criterion, one may iteratively double $n$ until the bounds are sufficiently tight. 

For lattices one tracks the decay of the complex-exponential Fourier coefficients \cite{JimHic16a}. For digital sequences one tracks the decay of  Fourier-Walsh coefficients  \cite{HicJim16a}. In both cases, the definitions of the node sets paired with the appropriate Fourier basis allows these approximate Fourier coefficients to be computed from the data, $\{y_i = f(\bsx_i)\}_{i=0}^{n -1} $ at a cost of $\calO(n \log n)$ for $n$ a power of $2$ using fast  transforms. 

Another idea which avoids multiple randomizations is to treat $f$ as a Gaussian process, as mentioned in Section \ref{sec:randBayes}. Conditioned on observations $f(\bsx_i)$, a closed form expression exists for the $1-\alpha$ credible intervals for $\mu$. While also applicable when using IID nodes, using LD sequences paired with special kernels enables these credible intervals to be computed at only $\calO(n \log n)$ cost \cite{Jag19a,JagHic22a,RatHic19a}. Again one may iteratively double the sample size $n$ until the credible intervals are sufficiently tight.

For this Bayesian approach, one must find the eigenvalues and eigenvectors of the Gram matrix $\mathsf{K} : = \bigl(K(\bsx_i,\bsx_j) \bigr)_{i,j=0}^{n-1}$.  In general this requires $\mathcal{O}(n^3)$ operations.  However, pairing shifted lattices with shift-invariant kernels as defined in \eqref{eq:shinvK} yields circulant Gram matrices for which computations can be done quickly using the fast Fourier transform \cite{RatHic19a}. Similarly, pairing digitally shifted digital nets with digitally-shift-invariant kernels yields recursive-block Toeplitz Gram matrices for which computations can be done quickly using the fast Walsh-Hadamard transform \cite{JagHic22a}. 

Extensions of these stopping criteria to functions of multiple expectations was explored in \cite{sorokin2022bounding}. Such problems include computing the posterior mean in a Bayesian context or computing  sensitivity indices which  quantifying feature importance to variability in a simulation. This work also considered sharing samples across approximations and generalized error criteria, such as a relative error tolerance.  

\section{Reformulating Our Problem} \label{sec:reformulate}

There are two themes of this section.  One is that our original problem may not look like $\mu = f(\bsX)$ for $\bsX \sim \mathcal{U}[0,1)^d$, and thus it may need to be rewritten.  The second is that qMC does best when most of the variation of $f$ is concentrated in the lower numbered coordinates, and, if possible, we should design $f$ so that this is so.

\subsection{Variable Transformations} 

Often qMC problems may be written as $\mu = \bbE[g(\bsT)]$ where $\bsT$ is a $d$-dimensional random variable with density $\lambda$. Since qMC nodes have LD with the standard uniform distribution $\bsX \sim \calU[0,1]^d$, we perform an invertible variable transformation $\bsPsi: [0,1]^d \to \calM$ satisfying $\mathrm{range}(\bsT) \subseteq \calM$. This implies
$\mu = \bbE[f(\bsX)]$ where  
$$
f(\bsX) = g(\bsPsi(\bsX)) \frac{\lambda (\bsPsi(\bsX))}{\nu(\bsPsi(\bsX))},
$$
and $\nu$ is the density of $\bsPsi(\bsX)$. Specifically, $\nu(\bsPsi(\bsX))= \lvert \bsJ_\bsPsi(\bsX) \rvert^{-1}$ where $\lvert \bsJ_\bsPsi(\bsX) \rvert$ is the determinant of the Jacobian of $\bsPsi$. 

The above formulation is equivalent to importance sampling by a distribution with density $\nu$. Optimally, one would choose $\nu \propto g\lambda$ so that $f$ is constant, but this is typically infeasible in practice due to our limited knowledge of $g$. Choosing $\bsPsi$ so that $\bsT = \bsPsi(\bsX)$ gives the simplification $f(\bsX) = g(\bsPsi(\bsX))$. An example of this was given in Section \ref{sec:practice} where $\bsT \sim \calN(\boldsymbol{0},\mathsf{I}/2)$ and $\bsPsi(\bsX) = \Phi^{-1}(\bsX)/\sqrt{2}$ where $\Phi^{-1}$ is the inverse CDF of the standard Gaussian distribution applied elementwise.

As observed in Section \ref{sec:coordwts}, qMC is more effective if the variation of $f$ is concentrated in the lower indexed coordinates.  For finance examples where the outcome depends more on the gross behavior of the underlying stochastic process, rather than the high frequency behavior, it makes sense to choose the transformation accordingly.  For example, if $\bsT$ represents a Brownian motion, then a principal component analysis coupled with the inverse Gaussian distribution shows greater convergence gains compared to the Cholesky decomposition.

\subsection{Control Variates} 

Since the upper bound of the qMC error includes the variation of $f$ as a factor, one may improve the efficiency of qMC if one can reformulate the problem to reduce this variation.  Choosing a good variable transformation, as mentioned in the previous subsection, is one approach.  

Another approach is control variates.  If $f_{\text{ctrl}} : [0,1)^d \to \RR$ is a function whose integral is zero, then $\mu(f) = \mu(f - f_{\text{ctrl}})$, and one may approximate $\mu(f)$ by $\hat{\mu}_n(f - f_{\text{ctrl}})$.  If the variation of  $f - f_{\text{ctrl}}$ is smaller than that of $f$, then it is likely that $\hat{\mu}_n(f - f_{\text{ctrl}})$ will have smaller error than $\hat{\mu}_n(f)$.  

An example of this in quantitative finance is where one may use the payoff of a European option, whose expectation is known, as a control variate for more exotic options.  As noted in \cite{HicEtal03}, reducing the variation of $f$ via control variates may or may not resemble reducing the variance of $f$.  Control variates have been implemented in the stopping criteria based on Fourier coefficients described in the previous section \cite{HicEtal17a}.

\subsection{Multilevel Methods} \label{sec:multi}
The cost of computing the sample mean defined in \eqref{eq:sample_mean} is typically $\mathcal{O}(dn)$ as the dimension, $d$, and the sample size, $n$, tend to infinity.  Some applications have a $d$ that is either huge or infinite.  On example is in quantitative finance, where $d$ denotes the number of times that the asset price is monitored over the life of a financial derivative times the number of assets involved.  

However, the $\mathcal{O}(dn)$ cost of approximating $\mu$ may often be reduced by employing multilevel methods.  In fact, these methods may reduce the cost to the level comparable to solving a low dimensional approximation \cite{Gil15a}. We introduce these methods here.

Let $\{Y_l := f_l(\bsX^{(l)})\}_{l=1}^L$ be a sequence of approximations to $Y = f(\bsX)$, where $\bsX^{(l)} = (X_1, \ldots, X_l) \sim \mathcal{U}[0,1)^{d_l}$.  Here we assume that $d_1 < \cdots < d_L$, and that the cost of evaluating $f_l(\bsx^{(l)})$ is $\mathcal{O}(d_l)$, where $\bsx^{(l)} = (x_1, \ldots, x_l)$.  Also, $d_L = d$, $Y_{d_L} = Y$, and $f_L = f$. Then the population mean may written as a telescoping sum:
\begin{align*}
    \mu & := \mathbb{E}(Y) = \mathbb{E}\bigl[f(\bsX)\bigr]
    = \mathbb{E}(Y_1) + \mathbb{E}(Y_2 - Y_1) + \cdots + \mathbb{E}(Y_{d_L} - Y_{d_{L-1}}) \\
    & = \int_{[0,1)^{d_1}} f_1(\bsx^{(1)}) \, \mathrm{d} \bsx^{(1)} +  \int_{[0,1)^{d_2}} \bigl[f_2(\bsx^{(2)}) - f_1(\bsx^{(1)})\bigr] \, \mathrm{d} \bsx^{(2)} + \cdots \\
    & \qquad \qquad + \int_{[0,1)^{d_L}} \bigl[f_L(\bsx^{(L)}) - f_{L-1}(\bsx^{(L-1)}) \bigr] \, \mathrm{d} \bsx^{(L)}.
\end{align*}
This series of integrals may be approximated by a series of sample means, namely
\begin{multline} \label{eq:MLmuhat}
    \hat{\mu}^{\text{ML}}_{\bsn,\bsd}  := \frac{1}{n_1}  \sum_{i=0}^{n_1-1} f_1(\bsx^{(1)}_i)  +  \frac{1}{n_2}  \sum_{i=0}^{n_2-1}  \bigl[f_2(\bsx^{(2)}_i) - f_1(\bsx^{(1)}_i)\bigr] + \cdots \\
     + \frac{1}{n_L}  \sum_{i=0}^{n_L-1}  \bigl[f_L(\bsx_i^{(L)}) - f_{L-1}(\bsx_i^{(L-1)}) \bigr].
\end{multline}
If $C_1$ is the computational cost of $f^{(1)}(\bsx^{(1)})$ and $C_l$ is the computational cost of $f^{(l)}(\bsx^{(1)}) - f_{l-1}(\bsx_i^{(l-1)})$ for $l = 2, \ldots, L$, then the total computational cost of 
\eqref{eq:MLmuhat} is 
\begin{equation} \label{eq:MLcostdef}
    \text{cost}(\bsn,\bsd) = n_1 C_1 + \cdots + n_L C_L.
\end{equation}

A worst case error analysis as described in Section \ref{sec:discrepancy} implies that
\begin{align}
    \label{eq:MLerrbd}
    \bigl \lvert \mu - \hat{\mu}^{\text{ML}}_{\bsn,\bsd} \bigr \rvert &\le
    V_1/n_1 + V_2/n_2 + \cdots + V_L/n_L,  \\
    \intertext{where the $V_l$ are defined such that we expect them to have little dependence on $n_l$ since we are using quasi-Monte Carlo nodes:} 
    \nonumber
    V_1 &= n_1 \, \text{discrepancy}\bigl(\{\bsx_i^{(1)}\}_{i=0}^{n_2-1}, K^{(1)} \bigr) \, \text{variation}(f^{(1)}, K^{(1)})\\
    \nonumber
    V_l & = n_l \,\text{discrepancy}\bigl(\{\bsx_i^{(l)}\}_{i=0}^{n_l-1}, K^{(l)} \bigr) \, \text{variation}(f^{(l)} - f^{(l-1)}, K^{(l)}), \\
    \nonumber
    & \qquad \qquad l = 2, \ldots, L,
\end{align}
and $K^{(l)}$ is the reproducing kernel corresponding to the relevant  Hilbert space of functions of $d_l$ variables.  

Optimizing $\bsn$ to minimize the error bound in \eqref{eq:MLerrbd} for a fixed cost and assuming constant $\bsC$ and $\bsV$, yields
\begin{align}
    \bigl \lvert \mu - \hat{\mu}^{\text{ML}}_{\bsn,\bsd} \bigr \rvert &\le \varepsilon \\
\intertext{for}
    n_l & = \frac{\sqrt{V_l} \bigl(\sqrt{V_1 C_1} + \ldots + \sqrt{V_LC_L}\bigr)}{\varepsilon \sqrt{C_l} }, \qquad l = 1, \ldots, d,\\
    \text{cost}(\bsn,\bsd) & = \frac{\bigl(\sqrt{V_1 C_1} + \ldots + \sqrt{V_LC_L}\bigr)^2}{\varepsilon}. \label{eq:MLcost}
\end{align}

We see that the cost of approximating $\mu$ within an absolute error tolerance of $\varepsilon$ is $\mathcal{O}(\varepsilon^{-1})$, as expected, but the dependence of $\text{cost}(\bsn,\bsd)$ on the dimension of the problem depends critically on how quickly the variation at each level decays.  Suppose that 
\begin{itemize}
    \item The cost of a function value at each level, $l$, is proportional to the dimension, $d_l$, and that the dimensions \emph{increase} exponentially, i.e.,  $d_l = d_1^{l}$ and $C_l \le \alpha d_1^l$. 
    \item The $V_l$, which are roughly the variation, \emph{decrease} exponentially, i.e., $V_l \le \beta r^{-l} $.
\end{itemize}   
Then the computational cost in \eqref{eq:MLcost} can be bounded above by
\begin{align*}
     \text{cost}(\bsn,\bsd) & \le \frac{\alpha \beta \bigl((d_1/r)^{1/2} + \ldots + (d_1/r)^{L/2}\bigr)^2}{\varepsilon} \\
     & = \frac{\alpha \beta \bigl((d_1/r)^{1/2} - (d_1/r)^{(L+1)/2}\bigr)^2}{\varepsilon \bigl( 1 - (d_1/r)^{1/2}\bigr)^2}.
\end{align*}

If the $V_l$ decay quickly enough, which means that the corrections $f_l - f_{l-1}$ depending on many variables decay quickly enough, in particular, if $r > d_1$, then 
\begin{equation}
         \text{cost}(\bsn,\bsd) \le \frac{\alpha \beta d_1} {r \varepsilon \bigl( 1 - (d_1/r)^{1/2}\bigr)^2} = \mathcal{O}(d_1/\varepsilon).
\end{equation}
In this case $\text{cost}(\bsn,\bsd)$ for this high dimensional problem is not much worse than the cost of ensuring that just error of the rough $d_1$ approximation is smaller than $\varepsilon$, which would be $V_1C_1/\varepsilon \le \alpha \beta d_1/(r\varepsilon)$.  

However, if $r < d_1$, then $\text{cost}(\bsn,\bsd) \ge V_LC_L/\varepsilon$ which may be as large as $\alpha \beta (d_1/r)^{-L}{\varepsilon}$, which grows exponentially as $d = d_L = d_1^L$ increases via increasing $L$.  Granted, this may not be as high a cost as the standard single-level approach where the cost would be roughly $V_1C_L$ if the $V_l$ are decaying.

It is an art to write one's problem to have the favorable case of high or even infinite nominal dimension, $d$, be solvable with a computational cost comparable to the low dimensional case.  However, the multilevel method takes advantage of the principle noted in Section \ref{sec:coordwts}, namely that qMC thrives when problem have a low effective dimension, even if the nominal dimension is large.  

\section{Beyond Computing Expectations} \label{sec:beyond}

Expectations of the form $\mu \coloneqq \mathbb{E}(f(\mathbf{X}))$ appear in several prominent areas of mathematics and qMC methods have long been applied in these contexts. Namely, in mathematical finance many problems can be formulated as high-dimensional integrals, where the large number of dimensions arises from small time steps in time discretization and/or a large number of variables. In this setting, Monte Carlo methods are applied to randomly simulate the evolution of stock prices over time, and the expected value is then estimated by averaging over all sample paths. The qMC approach replaces the randomly generated paths by carefully constructed paths in an effort to reduce the computational burden on the required number of paths. Excellent survey works for qMC and randomized qMC applied to finance include \cite{CafMor96,LEc09,Lem04a,Jae02,wangsloan05}. 

Another fruitful domain for qMC methods is in partial differential equations (PDEs), used for computing the expectations of nonlinear functionals of solutions of certain classes of PDEs. This includes problems in areas like fluid dynamics or material science, where random coefficients in the PDE model uncertain properties; see \cite{KuoNuy16a,KuoEtal12a,GraEtal11a} and references therein for qMC applied to PDEs.

QMC methods have also found recent applications in scientific machine learning of PDE solutions. It has been shown that LD nodes are more efficient than IID nodes for training neural networks modeling PDE solution \cite{longo21,mishra21}. LD nodes have also been used to more efficiently fit Gaussian process regression models for operator learning the solution process of PDEs with uncertain coefficients  \cite{sorokin2024computationally}. There the computational speedup follows from pairing LD nodes with special kernels to induce nicely structured gram matrices, as used in the fast Bayesian cubature stopping criteria discussed in Section \ref{sec:stop}. 

Besides computing expectations, qMC methods have found many applications in other areas of the computational sciences. We lightly review these applications to highlight their wider practical impact and to motivate the further development of qMC theory and methodology.

For example, when considering a variable $Y$, it is often useful to determine its distribution, $F$, or the density, $\varrho$, i.e., 
\[
F(y) \coloneqq \mathbb{P}(Y\leq y)=\mathbb{P}(f(\mathbf{X})\leq y)=\int_{[0,1]^d} \mathbf{1}_{(-\infty,y]}(f(\mathbf{x}))d\mathbf{x}, \quad \varrho(y)\coloneqq F'(y).
\]
Here, $F(y)$ gives the cumulative probability that $Y$ does not exceed $y$, while $\varrho(y)$ represents the corresponding density function. Several recent works for density estimation use qMC \cite{LEcPuc20a,LEcuyer2022b,GilKuoSlo23a} alongside randomized qMC \cite{AbdEtal21a} and often incorporate kernel-based methods or rely on conditioning strategies (sometimes called preintegration).

Keller \cite{kell95} was the first to adopt qMC for rendering photorealistic images in computer graphics where the problem of simulating light transport can be mathematically represented as an integral over all possible light paths from source to camera. Many algorithms involving sampling with LD nodes have since been developed much further in \cite{kellprem12,shired08,Keller2013a,kolkel00,kell96,kell95}. (Remarkably, Monte Carlo methods have even been awarded with an Academy Award for visual effects in movie making \cite{veachthesis,oscars}.)

Path planning is a critical aspect of efficient robotic operation, defining of a route from an initial position to a target destination while considering an optimality criterion. Most existing approaches utilize LD nodes, such as the Halton sequence \cite{veldel14,zhong24} with more recent efforts \cite{chahine24} demonstrating that employing state-of-the-art LD nodes from \cite{ruschkirk24} can significantly enhance the efficiency and flexibility of a robotic arm, demonstrated in both toy and real-world tasks.

In summary, qMC methods have demonstrated their applicability and impact well beyond computing expectations, finding fruitful applications across fields such as computer graphics, robotics, and machine learning.

\section{Conclusion} \label{sec:conclusion}
Quasi-Monte Carlo methods, based on low discrepancy sequences, can approximate the mean of a function, $f$, of $d$ variables to within an error tolerance of $\varepsilon$ at a cost of roughly $\mathcal{O}(\varepsilon^{-1})$ function evaluations, provided that the problem is formulated so that $f$ depends primarily on a modest number of variables.  This means that qMC outperforms alternative methods such as simple Monte Carlo and sampling on grids.  QMC is particularly powerful when the effective dimension of the problem is small, even though the nominal dimension, $d$, may be large.  Software libraries are available to facilitate qMC computations.

\bibliographystyle{spmpsci}
\bibliography{FJH25,FJHown23,main}

\def\Ignore#1{}\def\notesupp#1{}\providecommand{\HickernellFJ}{Hickernell\xspace}\def\Ignore#1{}\def\notesupp#1{}\providecommand{\HickernellFJ}{Hickernell\xspace}
\begin{thebibliography}{10}
\providecommand{\url}[1]{{#1}}
\providecommand{\urlprefix}{URL }
\expandafter\ifx\csname urlstyle\endcsname\relax
  \providecommand{\doi}[1]{DOI~\discretionary{}{}{}#1}\else
  \providecommand{\doi}{DOI~\discretionary{}{}{}\begingroup
  \urlstyle{rm}\Url}\fi

\bibitem{Aro50}
Aronszajn, N.: Theory of reproducing kernels.
\newblock Trans.\ Amer.\ Math.\ Soc. \textbf{68}, 337--404 (1950)

\bibitem{Atan04}
Atanassov, E.I.: On the discrepancy of the {H}alton sequences.
\newblock Math. Balkanica, New Series \textbf{18.1-2}, 15--32 (2004)

\bibitem{AbdEtal21a}
{Ben Abdellah}, A., L'Ecuyer, P., Owen, A.B., Puchhammer, F.: Density
  estimation by randomized quasi-{M}onte {C}arlo.
\newblock SIAM/ASA J.\ Uncertain.\ Quantif. \textbf{9}, 280--301 (2021).
\newblock \doi{10.1137/19M1259213}

\bibitem{BerT-A04}
Berlinet, A., Thomas-{A}gnan, C.: Reproducing Kernel {H}ilbert Spaces in
  Probability and Statistics.
\newblock Kluwer Academic Publishers, Boston (2004)

\bibitem{BriEtal18a}
Briol, F.X., Oates, C.J., Girolami, M., Osborne, M.A., Sejdinovic, D.:
  Probabilistic integration: A role in statistical computation?
\newblock Statist.\ Sci. \textbf{34}, 1--22 (2019)

\bibitem{CafMor96}
Caflisch, R.E., Morokoff, W.: Valuation of mortgage backed securities using the
  quasi-{M}onte {C}arlo method.
\newblock In: International Association of Financial Engineers First Annual
  Computational Finance Conference (1996)

\bibitem{chahine24}
Chahine, M., Rusch, T.K., Patterson, Z.J., Rus, D.: Improving efficiency of
  sampling-based motion planning via {M}essage-{P}assing {M}onte {C}arlo
  (2024).
\newblock \urlprefix\url{https://arxiv.org/abs/2410.03909}

\bibitem{MasChiWar05}
Chi, H., Mascagni, M., Warnock, T.: On the optimal {H}alton sequence.
\newblock Math. Comput. Simul. \textbf{70}(1), 9--21 (2005)

\bibitem{QMCPy2020a}
Choi, S.C.T., \HickernellFJ, F.J., Jagadeeswaran, R., McCourt, M., Sorokin, A.:
  {QMCPy}: A quasi-{M}onte {C}arlo {P}ython library (versions 1--1.5) (2024).
\newblock \doi{10.5281/zenodo.3964489}.
\newblock \urlprefix\url{https://qmcsoftware.github.io/QMCSoftware/}

\bibitem{clethesis24}
Cl{\'e}ment, F.: An optimization perspective on the construction of
  low-discrepancy point sets.
\newblock Ph.D. thesis, Sorbonne Universit\'{e} (2024)

\bibitem{cle24}
Cl{\'e}ment, F., Doerr, C., Klamroth, K., Paquete, L.: Constructing optimal
  {L}$_{\infty}$ star discrepancy sets (2024).
\newblock \urlprefix\url{https://arxiv.org/abs/2311.17463}

\bibitem{cle22}
Cl{\'e}ment, F., Doerr, C., Paquete, L.: Star discrepancy subset selection:
  Problem formulation and efficient approaches for low dimensions.
\newblock Journal of Complexity \textbf{70}, 101645 (2022)

\bibitem{cle24_heuristic}
Cl{\'e}ment, F., Doerr, C., Paquete, L.: Heuristic approaches to obtain
  low-discrepancy point sets via subset selection.
\newblock Journal of Complexity \textbf{83}, 101852 (2024)

\bibitem{CooNuy16a}
Cools, R., Nuyens, D. (eds.): {M}onte {C}arlo and Quasi-{M}onte {C}arlo
  Methods: {MCQMC}, {L}euven, {B}elgium, {A}pril 2014, \emph{Springer
  Proceedings in Mathematics and Statistics}, vol. 163. Springer-Verlag, Berlin
  (2016)

\bibitem{LatNet}
Darmon, Y., Godin, M., L'Ecuyer, P., Jemel, A., Marion, P., Munger, D.:
  {LatNet} builder (2018).
\newblock \urlprefix\url{https://github.com/umontreal-simul/latnetbuilder}

\bibitem{DicEtal22a}
Dick, J., Kritzer, P., Pillichshammer, F.: Lattice Rules: {N}umerical
  Integration, Approximation, and Discrepancy.
\newblock Springer Series in Computational Mathematics. Springer Cham (2022).
\newblock \doi{https://doi.org/10.1007/978-3-031-09951-9}

\bibitem{DicKuo04a}
Dick, J., Kuo, F.Y.: Reducing the construction cost of the
  component-by-component construction of good lattice rules.
\newblock Math.\ Comp. \textbf{73}, 1967---1988 (2004)

\bibitem{KuoEtal14a}
Dick, J., Kuo, F.Y., Peters, G.W., Sloan, I.H. (eds.): {M}onte {C}arlo and
  Quasi-{M}onte {C}arlo Methods 2012, \emph{Springer Proceedings in Mathematics
  and Statistics}, vol.~65. Springer-Verlag, Berlin (2013).
\newblock \doi{10.1007/978-3-642-41095-6}

\bibitem{DicPil10a}
Dick, J., Pillichshammer, F.: Digital Nets and Sequences: Discrepancy Theory
  and Quasi-{M}onte {C}arlo Integration.
\newblock Cambridge University Press, Cambridge (2010)

\bibitem{NieFanHic01}
Fang, K.T., \HickernellFJ, F.J., Niederreiter, H. (eds.): {M}onte {C}arlo and
  Quasi-{M}onte {C}arlo Methods 2000. Springer-Verlag, Berlin (2002)

\bibitem{faulem09}
Faure, H., Lemieux, C.: Generalized {H}alton sequences in 2008: A comparative
  study \textbf{19}(4) (2009)

\bibitem{GilKuoSlo23a}
Gilbert, A.D., Kuo, F.Y., Sloan, I.H.: Analysis of preintegration followed by
  quasi-{M}onte {C}arlo integration for distribution functions and densities.
\newblock SIAM J.\ Numer.\ Anal. \textbf{61}, 135--166 (2023)

\bibitem{Gil15a}
Giles, M.B.: Multilevel {M}onte {C}arlo methods.
\newblock Acta Numer. \textbf{24}, 259--328 (2015).
\newblock \doi{10.1017/S096249291500001X}

\bibitem{GraEtal11a}
Graham, I.G., Kuo, F.Y., Nuyens, D., Scheichl, R., Sloan, I.H.: Quasi-{M}onte
  {C}arlo methods for elliptic {PDEs} with random coefficients and
  applications.
\newblock J. Comput.\ Phys. \textbf{230}, 3668--3694 (2011).
\newblock \doi{doi.org/10.1016/j.jcp.2011.01.023}

\bibitem{HeiHicYue02a}
Heinrich, S., \HickernellFJ, F.J., Yue, R.X.: Optimal quadrature for {H}aar
  wavelet spaces.
\newblock Math.\ Comp. \textbf{73}, 259--277 (2004)

\bibitem{MCQMC2024TutorialNotebook}
Hickernell, F.J., Kirk, N., Sorokin, A.G.:  (2025).
\newblock
  \urlprefix\url{https://github.com/QMCSoftware/QMCSoftware/blob/MCQMC2024/demos/talk_paper_demos/MCQMC2024Tutorial/MCQMC2024Tutorial.ipynb}

\bibitem{Hic95}
\HickernellFJ, F.J.: A comparison of random and quasirandom points for
  multidimensional quadrature.
\newblock In: Niederreiter and Shiue  \cite{NieShi95}, pp. 213--227

\bibitem{Hic97a}
\HickernellFJ, F.J.: A generalized discrepancy and quadrature error bound.
\newblock Math.\ Comp. \textbf{67}, 299--322 (1998).
\newblock \doi{10.1090/S0025-5718-98-00894-1}

\bibitem{Hic98b}
\HickernellFJ, F.J.: Lattice rules: How well do they measure up?
\newblock In: P.~Hellekalek, G.~Larcher (eds.) Random and Quasi-Random Point
  Sets, \emph{Lecture Notes in Statistics}, vol. 138, pp. 109--166.
  Springer-Verlag, New York (1998)

\bibitem{Hic99b}
\HickernellFJ, F.J.: What affects the accuracy of quasi-{M}onte {C}arlo
  quadrature?
\newblock In: H.~Niederreiter, J.~Spanier (eds.) {M}onte {C}arlo and
  Quasi-{M}onte {C}arlo Methods 1998, pp. 16--55. Springer-Verlag, Berlin
  (2000)

\bibitem{Hic00a}
\HickernellFJ, F.J., Hong, H.S.: Quasi-{M}onte {C}arlo methods and their
  randomizations.
\newblock In: R.~Chan, Y.K. Kwok, D.~Yao, Q.~Zhang (eds.) Applied Probability,
  \emph{AMS/IP Studies in Advanced Mathematics}, vol.~26, pp. 59--77. American
  Mathematical Society, Providence, Rhode Island (2002)

\bibitem{HicEtal00}
\HickernellFJ, F.J., Hong, H.S., L'{\'E}cuyer, P., Lemieux, C.: Extensible
  lattice sequences for quasi-{M}onte {C}arlo quadrature.
\newblock SIAM J. Sci.\ Comput. \textbf{22}, 1117--1138 (2000).
\newblock \doi{10.1137/S1064827599356638}

\bibitem{HicEtal14a}
\HickernellFJ, F.J., Jiang, L., Liu, Y., Owen, A.B.: Guaranteed conservative
  fixed width confidence intervals via {M}onte {C}arlo sampling.
\newblock In: Dick et~al.  \cite{KuoEtal14a}, pp. 105--128.
\newblock \doi{10.1007/978-3-642-41095-6}

\bibitem{HicJim16a}
\HickernellFJ, F.J., {Jim\'enez Rugama}, {\relax Ll}.A.: Reliable adaptive
  cubature using digital sequences.
\newblock In: Cools and Nuyens  \cite{CooNuy16a}, pp. 367--383.
\newblock ArXiv:1410.8615 [math.NA]

\bibitem{HicEtal17a}
\HickernellFJ, F.J., {Jim\'enez Rugama}, {\relax Ll}.A., Li, D.: Adaptive
  quasi-{M}onte {C}arlo methods for cubature.
\newblock In: J.~Dick, F.Y. Kuo, H.~Wo\'zniakowski (eds.) Contemporary
  Computational Mathematics --- a celebration of the 80th birthday of {I}an
  {S}loan, pp. 597--619. Springer-Verlag (2018).
\newblock \doi{10.1007/978-3-319-72456-0}

\bibitem{HicEtal03}
\HickernellFJ, F.J., Lemieux, C., Owen, A.B.: Control variates for
  quasi-{M}onte {C}arlo.
\newblock Statist.\ Sci. \textbf{20}, 1--31 (2005).
\newblock \doi{10.1214/088342304000000468}

\bibitem{HicNie03a}
\HickernellFJ, F.J., Niederreiter, H.: The existence of good extensible rank-1
  lattices.
\newblock J. Complexity \textbf{19}, 286--300 (2003)

\bibitem{HicYue00}
\HickernellFJ, F.J., Yue, R.X.: The mean square discrepancy of scrambled
  $(t,s)$-sequences.
\newblock SIAM J. Numer.\ Anal. \textbf{38}, 1089--1112 (2000).
\newblock \doi{10.1137/S0036142999358019}

\bibitem{hinoet16}
Hinrichs, A., Oettershagen, J.: Optimal point sets for quasi-{M}onte {C}arlo
  integration of bivariate periodic functions with bounded mixed derivatives.
\newblock In: Monte Carlo and Quasi-Monte Carlo Methods 2016, pp. 385--405
  (2018)

\bibitem{Jae02}
J{\"a}ckel, P.: Monte {C}arlo Methods in Finance.
\newblock John Wiley \& Sons Ltd., Chichester, England (2002)

\bibitem{Jag19a}
Jagadeeswaran, R.: Fast automatic bayesian cubature using matching kernels and
  designs.
\newblock Ph.D. thesis, Illinois Institute of Technology (2019)

\bibitem{RatHic19a}
Jagadeeswaran, R., \HickernellFJ, F.J.: Fast automatic {B}ayesian cubature
  using lattice sampling.
\newblock Stat.\ Comput. \textbf{29}, 1215--1229 (2019).
\newblock \doi{10.1007/s11222-019-09895-9}

\bibitem{JagHic22a}
Jagadeeswaran, R., \HickernellFJ, F.J.: Fast automatic {B}ayesian cubature
  using {S}obol' sampling.
\newblock In: Z.~Botev, A.~Keller, C.~Lemieux, B.~Tuffin (eds.) Advances in
  Modeling and Simulation: {F}estschrift in Honour of {P}ierre {L'E}cuyer, pp.
  301--318. Springer, Cham (2022).
\newblock \doi{10.1007/978-3-031-10193-9\_15}

\bibitem{JimHic16a}
{Jim\'enez Rugama}, {\relax Ll}.A., \HickernellFJ, F.J.: Adaptive
  multidimensional integration based on rank-1 lattices.
\newblock In: Cools and Nuyens  \cite{CooNuy16a}, pp. 407--422.
\newblock ArXiv:1411.1966

\bibitem{JoeKuo03}
Joe, S., Kuo, F.Y.: Remark on algorithm 659: Implementing {S}obol's quasirandom
  sequence generator.
\newblock ACM Trans.\ Math.\ Software \textbf{29}, 49--57 (2003)

\bibitem{KuoJoe08a}
Joe, S., Kuo, F.Y.: Constructing {S}obol sequences with better two-dimensional
  projections.
\newblock SIAM J.\ Sci.\ Comput. \textbf{30}, 2635--2654 (2008)

\bibitem{Kei96}
Keister, B.D.: Multidimensional quadrature algorithms.
\newblock Computers in Physics \textbf{10}, 119--122 (1996).
\newblock \doi{10.1063/1.168565}

\bibitem{kell95}
Keller, A.: A quasi-{M}onte {C}arlo algorithm for the global illumination
  problem in a radiosity setting.
\newblock In: In Harald Niederreiter and Peter Jau-Shyong Shiue, editors, Monte
  Carlo and Quasi-Monte Carlo Methods in Scientific Computing, pp. 239--251.
  Springer-Verlag, New York (1995)

\bibitem{kell96}
Keller, A.: Quasi-{M}onte {C}arlo radiosity.
\newblock In: X. Pueyo and P. Schr{\"o}der (Eds.), Eurographics Rendering
  Workshop, pp. 101--110 (1996)

\bibitem{Keller2013a}
Keller, A.: Quasi-{M}onte {C}arlo image synthesis in a nutshell.
\newblock In: Dick et~al.  \cite{KuoEtal14a}, pp. 213--249.
\newblock \doi{10.1007/978-3-642-41095-6}

\bibitem{Kel22a}
Keller, A. (ed.): {M}onte {C}arlo and Quasi-{M}onte {C}arlo Methods: {MCQMC},
  {O}xford, England, {A}ugust 2020, Springer Proceedings in Mathematics and
  Statistics. Springer, Cham (2022)

\bibitem{kellprem12}
Keller, A., Premoze, S., Raab, M.: Advanced (quasi) {M}onte {C}arlo methods for
  image synthesis.
\newblock In: ACM SIGGRAPH 2012 Courses, SIGGRAPH '12. Association for
  Computing Machinery, New York, NY, USA (2012)

\bibitem{kirklem24}
Kirk, N., Lemieux, C.: An improved {H}alton sequence for implementation in
  quasi-{M}onte {C}arlo methods.
\newblock In: 2023 Winter Simulation Conference (WSC). IEEE (2024)

\bibitem{kolkel00}
Kollig, T., Keller, A.: Efficient bidirectional path tracing by randomized
  quasi-{M}onte {C}arlo integration.
\newblock In: Fang et~al.  \cite{NieFanHic01}, pp. 290--305

\bibitem{kor63}
Korobov, N.M.: Number-theoretic methods in approximate analysis.
\newblock Goz. Izdat. Fiz.-Math.  (1963)

\bibitem{KuoJoe02b}
Kuo, F.Y.: Component-by-component constructions achieve the optimal rate of
  convergence for multivariate integration in weighted {K}orobov and {S}obolev
  spaces.
\newblock J.\ Complexity \textbf{19}, 301--320 (2003)

\bibitem{KuoNuy16a}
Kuo, F.Y., Nuyens, D.: Application of quasi-{M}onte {C}arlo methods to elliptic
  {PDE}s with random diffusion coefficients -- a survey of analysis and
  implementation.
\newblock Found.\ Comput.\ Math. \textbf{16}, 1631--1696 (2016).
\newblock \urlprefix\url{https://people.cs.kuleuven.be/~dirk.nuyens/qmc4pde/}

\bibitem{KuoEtal12a}
Kuo, F.Y., Schwab, C., Sloan, I.H.: Quasi-{M}onte {C}arlo finite element
  methods for a class of elliptic partial differential equations with random
  coefficients.
\newblock SIAM J.\ Numer.\ Anal. \textbf{50}, 3351--3374 (2012).
\newblock \doi{10.1137/110845537}

\bibitem{LEc09}
L'Ecuyer, P.: Quasi-{M}onte {C}arlo methods with applications in finance.
\newblock Finance Stoch. \textbf{13}, 307--349 (2009)

\bibitem{LEcEtal22a}
L'Ecuyer, P., Marion, P., Godin, M., Puchhammer, F.: A tool for custom
  construction of {QMC} and {RQMC} point sets.
\newblock In: Keller  \cite{Kel22a}

\bibitem{LEcPuc20a}
L'Ecuyer, P., Puchhammer, F.: Density estimation by {M}onte {C}arlo and
  quasi-{M}onte {C}arlo.
\newblock In: Keller  \cite{Kel22a}, pp. 3--21

\bibitem{LEcuyer2022b}
L'Ecuyer, P., Puchhammer, F., {Ben Abdellah}, A.: {M}onte {C}arlo and
  quasi{\textendash}{M}onte {C}arlo density estimation via conditioning.
\newblock INFORMS J.\ Comput. \textbf{34}(3), 1729--1748 (2022).
\newblock \doi{10.1287/ijoc.2021.1135}.
\newblock \urlprefix\url{https://doi.org/10.1287/ijoc.2021.1135}

\bibitem{Lem04a}
Lemieux, C.: Randomized quasi-{M}onte {C}arlo: a tool for improving the
  efficiency of simulations in finance.
\newblock In: R.G. Ingalls, M.D. Rossetti, J.S. Smith, B.A. Peters (eds.)
  Proc.\ 2004 Winter Simulation Conference, pp. 1565--1573. IEEE Press (2004)

\bibitem{longo21}
Longo, M., Mishra, S., Rusch, T.K., Schwab, C.: Higher-order quasi-{M}onte
  {C}arlo training of deep neural networks.
\newblock SIAM Journal on Scientific Computing \textbf{43}(6), A3938--A3966
  (2021)

\bibitem{l2023confidence}
L’Ecuyer, P., Nakayama, M.K., Owen, A.B., Tuffin, B.: Confidence intervals
  for randomized quasi-{M}onte {C}arlo estimators.
\newblock In: 2023 Winter Simulation Conference (WSC), pp. 445--456. IEEE
  (2023)

\bibitem{Mai81a}
Maize, E.: Contributions to the theory of error reduction in quasi-{M}onte
  {C}arlo methods.
\newblock Ph.D. thesis, The Claremont Graduate School (1981)

\bibitem{Mat98}
Matou\v{s}ek, J.: On the {$L_2$}-discrepancy for anchored boxes.
\newblock J.\ Complexity \textbf{14}, 527--556 (1998)

\bibitem{mishra21}
Mishra, S., Rusch, T.K.: Enhancing accuracy of deep learning algorithms by
  training with low-discrepancy sequences.
\newblock SIAM Journal on Numerical Analysis \textbf{59}(3), 1811--1834 (2021)

\bibitem{Nie92}
Niederreiter, H.: Random Number Generation and Quasi-{M}onte {C}arlo Methods.
\newblock CBMS-NSF Regional Conference Series in Applied Mathematics. SIAM,
  Philadelphia (1992)

\bibitem{NieShi95}
Niederreiter, H., Shiue, P.J.S. (eds.): {M}onte {C}arlo and Quasi-{M}onte
  {C}arlo Methods in Scientific Computing, \emph{Lecture Notes in Statistics},
  vol. 106. Springer-Verlag, New York (1995)

\bibitem{NovWoz10a}
Novak, E., Wo{\'{z}}niakowski, H.: Tractability of Multivariate Problems
  {V}olume {II}: {S}tandard Information for Functionals.
\newblock No.~12 in EMS Tracts in Mathematics. European Mathematical Society,
  Z\"urich (2010)

\bibitem{NuyCoo06b}
Nuyens, D., Cools, R.: Fast algorithms for component-by-component construction
  of rank-1 lattice rules in shift-invariant reproducing kernel {H}ilbert
  spaces.
\newblock Math.\ Comp. \textbf{75}, 903---920 (2006)

\bibitem{NuyCoo06a}
Nuyens, D., Cools, R.: Fast component-by-component construction.
\newblock In: H.~Niederreiter, D.~Talay (eds.) {M}onte {C}arlo and
  Quasi-{M}onte {C}arlo Methods 2004, pp. 373--387. Springer-Verlag, Berlin
  (2006)

\bibitem{oscars}
Oscars: The 86th scientific \& technical awards 2013 | 2014 (2014).
\newblock \urlprefix\url{https://www.oscars.org/sci-tech/ceremonies/2014}

\bibitem{Owe95}
Owen, A.B.: Randomly permuted $(t,m,s)$-nets and $(t,s)$-sequences.
\newblock In: Niederreiter and Shiue  \cite{NieShi95}, pp. 299--317

\bibitem{Owe97}
Owen, A.B.: Scrambled net variance for integrals of smooth functions.
\newblock Ann.\ Stat. \textbf{25}, 1541--1562 (1997)

\bibitem{owen2024gain}
Owen, A.B., Pan, Z.: Gain coefficients for scrambled {H}alton points.
\newblock SIAM Journal on Numerical Analysis \textbf{62}(3), 1021--1038 (2024)

\bibitem{PanOwe23a}
Pan, Z., Owen, A.B.: Super-polynomial accuracy of one dimensional randomized
  nets using the median-of-means.
\newblock Math.\ Comp. \textbf{92}, 805--837 (2023).
\newblock \doi{10.1090/mcom/3791}

\bibitem{paulin2022}
Paulin, L., Bonneel, N., Coeurjolly, D., Iehl, J.C., Keller, A., Ostromoukhov,
  V.: Matbuilder: Mastering sampling uniformity over projections.
\newblock ACM Transactions on Graphics (Proceedings of SIGGRAPH) \textbf{41}(4)
  (2022).
\newblock \doi{https://doi.org/10.1145/3528223.3530063}

\bibitem{Rit00a}
Ritter, K.: Average-Case Analysis of Numerical Problems, \emph{Lecture Notes in
  Mathematics}, vol. 1733.
\newblock Springer-Verlag, Berlin (2000)

\bibitem{ruschkirk24}
Rusch, T.K., Kirk, N., Bronstein, M.M., Lemieux, C., Rus, D.: Message-{P}assing
  {M}onte {C}arlo: Generating low-discrepancy point sets via graph neural
  networks.
\newblock Proceedings of the National Academy of Sciences \textbf{121}(40),
  e2409913121 (2024)

\bibitem{shired08}
Shirley, P., Edwards, D., Boulos, S.: {M}onte {C}arlo and quasi-{M}onte {C}arlo
  methods for computer graphics.
\newblock In: A.~Keller, S.~Heinrich, H.~Niederreiter (eds.) Monte Carlo and
  Quasi-Monte Carlo Methods 2006, pp. 167--177. Springer-Verlag, Berlin (2008)

\bibitem{Slo02}
Sloan, I.H.: {QMC} integration --- beating intractability by weighting the
  coordinate directions.
\newblock In: Fang et~al.  \cite{NieFanHic01}, pp. 103--123

\bibitem{SloJoe94}
Sloan, I.H., Joe, S.: Lattice Methods for Multiple Integration.
\newblock Oxford University Press, Oxford (1994)

\bibitem{SloWoz98}
Sloan, I.H., Wo\'zniakowski, H.: When are quasi-{M}onte {C}arlo algorithms
  efficient for high dimensional integrals?
\newblock J. Complexity \textbf{14}, 1--33 (1998)

\bibitem{Sob67}
Sobol', I.M.: The distribution of points in a cube and the approximate
  evaluation of integrals.
\newblock U.S.S.R. Comput. Math. and Math. Phys. \textbf{7}, 86--112 (1967)

\bibitem{sorokin2024computationally}
Sorokin, A.G., Pachalieva, A., O’Malley, D., Hyman, J.M., Hickernell, F.J.,
  Hengartner, N.W.: Computationally efficient and error aware surrogate
  construction for numerical solutions of subsurface flow through porous media.
\newblock Advances in Water Resources \textbf{193}, 104836 (2024)

\bibitem{sorokin2022bounding}
Sorokin, A.G., Rathinavel, J.: On bounding and approximating functions of
  multiple expectations using quasi-{M}onte {C}arlo.
\newblock In: International Conference on Monte Carlo and Quasi-Monte Carlo
  Methods in Scientific Computing, pp. 583--599. Springer (2022)

\bibitem{tuffin98}
Tuffin, B.: A new permutation choice in {H}alton sequences.
\newblock In: Proceedings of the Conference on Monte Carlo and Quasi-Monte
  Carlo Methods, P. Hellekalek, G. Larcher, H. Niederreiter, and P. Zinterhof,
  Eds. Lecture Notes in Statistics, vol. 127. Springer, Berlin, New York, pp.
  427--435 (1998)

\bibitem{vancools06}
Vandewoestyne, B., Cools, R.: Good permutations for deterministic scrambled
  {H}alton sequences in terms of {L}2-discrepancy.
\newblock Journal of Computational and Applied Mathematics \textbf{189},
  341--361 (2006).
\newblock Proceedings of The 11th International Congress on Computational and
  Applied Mathematics

\bibitem{veachthesis}
Veach, E.: Robust {M}onte {C}arlo methods for light transport simulation.
\newblock Ph.D. thesis, Stanford University (1997)

\bibitem{veldel14}
Velagi{\'c}, J., Delimustafi{\'c}, D., Osmankovi{\'c}, D.: Mobile robot
  navigation system based on probabilistic road map ({P}{R}{M}) with {H}alton
  sampling of configuration space.
\newblock In: 2014 IEEE 23rd International Symposium on Industrial Electronics
  (ISIE), pp. 1227--1232 (2014)

\bibitem{wangsloan05}
Wang, X., Sloan, I.H.: Why are high-dimensional finance problems often of low
  effective dimension?
\newblock SIAM J. Sci. Comput. \textbf{27 (1)}, 159--183 (2005)

\bibitem{WinFan97b}
Winker, P., Fang, K.T.: Optimal {U}-type designs.
\newblock In: H.~Niederreiter, P.~Hellekalek, G.~Larcher, P.~Zinterhof (eds.)
  {M}onte {C}arlo and quasi-{M}onte {C}arlo methods 1996, \emph{Lecture Notes
  in Statistics}, vol. 127, pp. 436--448. Springer-Verlag, New York (1998)

\bibitem{zhong24}
Zhong, H., Cong, M., Wang, M., Du, Y., Liu, D.: H{B}-{R}{R}{T}: {A} path
  planning algorithm for mobile robots using {H}alton sequence-based
  rapidly-exploring random tree.
\newblock Engineering Applications of Artificial Intelligence \textbf{133},
  108362 (2024)

\end{thebibliography}

\end{document}